\documentclass[11pt]{article}
\usepackage{lscape}
\usepackage[centertags]{amsmath}
\usepackage{amsfonts}
\usepackage{amssymb}
\usepackage{amsthm}
\usepackage{newlfont}
\newskip\ns \ns=1cm
\def\bkR{{\rm I\kern-.17em R}}

\date{}
\begin{document}
\title{On Characterizing weak defining hyperplanes (weak Facets) in DEA
with Constant Returns to Scale Technology}
\author{ Dariush Akbarian \footnote{Corresponding author,\newline{\it E-mail
address:} d\underline{ }akbarian@yahoo.com, d-akbarian@iau-arak.ac.ir {\it Fax:} 0098-861-3120532}}

 %%%----------------------------------------------------------------------
\maketitle
\begin{center}\textit{{\scriptsize  Department of Mathematics, Arak Branch, Islamic
Azad University, Arak, Iran.\\ }}\end{center}
%\begin{center}
%E-mail address: d\underline{ }akbarian@yahoo.com, d-akbarian@iau-arak.ac.ir\\
%Fax: 0098-861-3120532
%\end{center}
\begin{abstract}
The Production Possibility Set (PPS) is defined as a set of inputs and outputs of a system in which inputs can produce outputs. The Production Possibility Set of the Data Envelopment
Analysis (DEA) model is contain of two types defining hyperplanes (facets); strong and weak efficient facets. In this paper, the problem of finding weak defining hyperplanes of the PPS
of the CCR-technology is dealt with. We state and prove some properties relative to our method. To illustrate the applicability of the proposed model, some numerical examples are finally
provided. Our algorithm can easily be implemented using existing packages for operation research, such as GAMS.\\\\
%\textbf{\textit{AMS}}: \textrm
%{Operation Research, No. 90.}\\
\textbf{\textit{Keywords}}: \textrm{Data Envelopment Analysis (DEA); Production Possibility Set; Efficient frontier; Weak efficient frontier; Hyperplane.}
\end{abstract}

%%% -------------------------------------------------------------------
\section{Introduction}
The Data Envelopment Analysis (DEA), introduced by Charnes, Cooper and Rhodes (CCR) (1978), is a procedure to evaluate the relative efficiency of a set of decision making units (DMU); each uses multiple inputs to produce multiple outputs. The data define a Production Possibility Set (PPS); that can be used to evaluate the efficiency of each of DMUs. The PPS of the DEA models is the smallest set containing the observed DMUs and all feasible input–output level correspondences pertaining to the production process operated by the DMUs. The PPS of the CCR model is the intersection of a finite number of halfspace, whose defining hyperplanes pass through the origin. The defining hyperplanes (facets) of the PPS of the CCR model are divided into two categories, including {\it i}) strong defining hyperplanes, and {\it ii}) weak defining hyperplanes. One of the problems in DEA is to find the equations of these defining hyperplanes. There are many researches undertaken on the subject of finding {\it strong} defining hyperplanes (see for example Amirteimoori et al. (2005), Amirteimoori et al. (2012), M. Davtalab-Olyaie et. al. (2014), Jahanshahloo et al. (2007), Jahanshahloo et al. (2005a), Yu et al. (1996), Wei et al. (2007) and Olesen et al. (2003)). However, less attention has been paid about finding {\it weak} defining hyperplanes of the PPS of the CCR model (see Wei et al. (2007)). In this regard, Jahanshahloo et al. (2010), proposed a method for determining weak defining hyperplanes of the PPS of the BCC model. In this paper we provide a method to find {\it weak} defining hyperplanes of the PPS of the CCR model.

Using the method proposed by Jahanshahloo et al. (2007) (with some modifications) we introduce a method to find the equations of {\it weak} defining hyperplanes of the PPS of the CCR model. The idea to find weak defining hyperplanes is straightforward by adding artificial weak efficient DMUs, named weak efficient virtual DMUs in this paper. The key is how to define these artificial weak efficient DMUs, and it is done by testing all CCR-efficient DMUs by a variance of super-efficiency models (see models (\ref{tak}) and (\ref{tao}))(after eliminating all CCR-inefficient DMUs from the PPS) and determining all extreme DMUs that lie on the some weak efficient defining hyperplanes. A supporting hyperplane is found to be a weak defining hyperplane if at least one artificial DMU lies on it. Using this method, it is possible to check (i) which CCR-efficient DMUs lie on the extreme rays (edges) of the PPS of the CCR model, (ii) which extreme DMUs lie on the some weak defining hyperplanes of the PPS of the CCR model, (iii) how many weak and strong defining hyperplanes they are on. Moreover, there are many PPS of DEA models that don't have any strong defining hyperplanes (see Figure \ref{2} and remark 2). Also, these hyperplanes are useful in sensitivity and stability analysis (see Jahanshahloo et al. (2005b)). Finally, the alternative optimal solutions of the models (\ref{multii}) and (\ref{multio}) can be found using the proposed method. These may show the importance of obtaining the weak defining hyperplanes of the PPS of DEA model. Some useful facts related to the properties of models (\ref{tak}) and (\ref{tao}) are stated and proved. In addition, three numerical examples are provided.
\section{Background}
Consider a set of $\textbf{\textit{n}}$ DMUs which is associated with $\textbf{\textit{m}}$ inputs and $\textbf{\textit{s}}$ outputs. Particularly, each $DMU_{j}=(X_j, Y_j)$ $(j\in J=\{1,\ldots, n\})$ consumes amount $x_{ij}(>0)$ of input $\textbf{\textit{i}}$ and produces amount $y_{rj}(>0)$ of output $\textbf{\textit{r}}$. The production possibility set $T$, $T\subset\big\{(X, Y)|X\in E^{m}, Y\in E^{s}, X\geqslant0, Y\geqslant0\big\}$ is based on postulate sets which are presented with a brief explanation (see Banker (1984), Banker et al. (1984) and Yu et al. (1996)). One of the DEA models to evaluate the relative efficiency of a set of DMUs is the CCR model, which is, proposed by Charnes et al. (1978). The production possibility set (PPS) of the CCR model can be defined as follows:

$$T=\Big\{(X, Y)|X\geqq\displaystyle\sum_{j\in J}\lambda_{j}X_{j},\ Y\leqq\displaystyle\sum_{j\in J}\lambda_{j}Y_{j},\ \lambda_{j}\geqslant0,\ j\in J\Big\}.$$ in which $X_{j}$ and $Y_{j}$ are vectors of input and output of $DMU_{j}$, respectively.\\

Following properties are postulated for PPS \textit{T}:
\begin{description}
  \item[1.] The observed activities belongs to \textit{T}; i.e.\\ $(X_{j}, Y_{j},)\in T$, $j=1, ..., n$
 \item[2.] If $(X, Y)\in T$, then the $(tX, tY)\in T$ for any $t>0$.
  \item[3.] For any activity $(X, Y)\in T$ any nonnegative activity $(\bar{X}, \bar{Y})$ with $\bar{X}\geq X$ and $\bar{Y}\leq Y$ is included in \textit{T}.
  \item[4.] \textit{T} is closed and convex.
\end{description}
A \textit{face} of a polyhedral set is the support set of a supporting hyperplane.\\ A \textit{facet} of a $k$-dimensional polyhedral set is a $k-1$ dimensional face. In fact, any facet of the PPS of the DEA model is a defining hyperplane of the PPS.\\
\textit{Note}: A CCR-efficient DMU is said to be \textit{extreme} DMU; if it lies on the edge of the PPS of the CCR model.

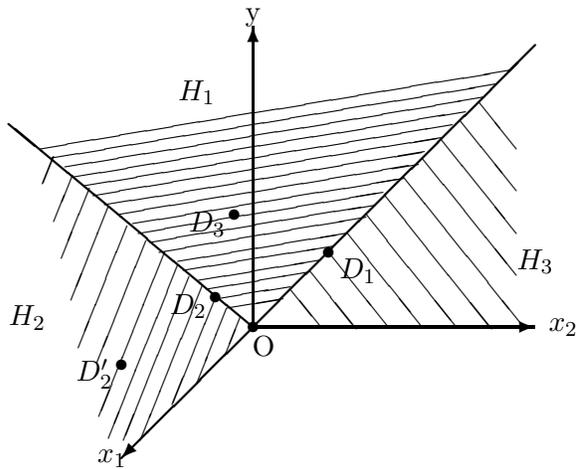
\begin{figure}
[t]
\begin{center}
\unitlength=2.5mm
\begin{picture}(18,18)(-4,-4)
\thicklines\put(0,0){\vector(1,0){15}}
\put(0,0){\vector(0,1){16}}
\thicklines\put(0,0){\vector(-1,-1){7}}
\put(0,0){\makebox(33,0)[]{\shortstack[b]{$x_{2}$}}}
\put(0,0){\makebox(0,33)[]{\shortstack[b]{y}}}
\put(0,0){\makebox(-15,-14)[]{\shortstack[b]{$x_{1}$}}}
\put(-1,6){{\circle*{.6}}\makebox(-4,-1)[]{\shortstack[b]{$D_{3}$}}}
\put(-1,13){\makebox(-4,-1)[]{\shortstack[b]{$H_1$}}}
\put(-7.,-2){{\circle*{.6}}\makebox(-4,-1)[]{\shortstack[b]{$D'_{2}$}}}
\put(-10,1){\makebox(-4,-1)[]{\shortstack[b]{$H_2$}}}
\put(-2,1.6){{\circle*{.6}}\makebox(-4,-1)[]{\shortstack[b]{$D_{2}$}}}
%\put(-2,-2){{\circle*{.6}}\makebox(-1,-2.4)[]{\shortstack[b]{$D'_{2}$}}}
\put(4,4){{\circle*{.6}}\makebox(2,-2)[]{\shortstack[b]{$D_{1}$}}}
\put(17,4){\makebox(-4,-1)[]{\shortstack[b]{$H_3$}}}
%\put(5,0){{\circle*{.6}}\makebox(2.5,-2)[]{\shortstack[b]{$D'_{1}$}}}
\put(0,0){{\circle*{.6}}\makebox(0,-2)[]{\shortstack[b]{O}}}
\thicklines\put(0,0){\line(1,1){15}}
\thicklines\put(0,0){\line(-6,5){13}}
\thinlines\put(-1.4,1){\line(6,1){3.}}\put(-2.,1.5){\line(6,1){4}}\put(-2.5,2){\line(6,1){5.5}}\put(-2.9,2.5){\line(6,1){6.5}}\put(-4,3.5){\line(6,1){9}}\put(-5,4){\line(6,1){11}}
\thinlines\put(-3.7,3){\line(6,1){8}}\put(-5.5,4.5){\line(6,1){12}}
\put(-6.1,5){\line(6,1){13}}\put(-6.4,5.5){\line(6,1){14.4}}
\thinlines\put(-7.2,6.){\line(6,1){16}}\put(-7.9,6.5){\line(6,1){17}}
\put(-8.5,7){\line(6,1){18.5}}\put(-9,7.5){\line(6,1){19.6}}
\put(-9.5,8){\line(6,1){20.8}}\put(-10,8.5){\line(6,1){22.3}}\put(-10.7,9){\line(6,1){23.6}}\put(-11.3,9.5){\line(6,1){25}}\put(-.7,.5){\line(-2,-5){.8}}\put(-1.4,1.){\line(-2,-5){1.6}}\put(-2,1.6){\line(-2,-5){2.3}}\put(-2.7,2.2){\line(-2,-5){3.3}}
\put(-3.4,2.7){\line(-2,-5){3.5}}\put(-4,3.3){\line(-2,-5){3.7}}\put(-4.7,4.){\line(-2,-5){3.5}}\put(-5.7,4.6){\line(-2,-5){2.6}}\put(-6.7,5.4){\line(-2,-5){2.6}}\put(-7.7,6.4){\line(-2,-5){2.}}
\put(-8.7,7.4){\line(-2,-5){1.5}}\put(-9.6,8.1){\line(-2,-5){1}}\put(-10.6,9.1){\line(-2,-5){.6}}
\put(1,1.2){\line(4,-5){.94}}\put(2,2){\line(4,-5){1.6}}\put(3,3){\line(4,-5){2.5}}\put(4,4){\line(4,-5){3.2}}\put(5,5){\line(4,-5){3.9}}
\put(6,6){\line(4,-5){4.8}}\put(7,7){\line(4,-5){5.7}}\put(8,8){\line(4,-5){6.4}}\put(9,9){\line(4,-5){5}}\put(10,10){\line(4,-5){4}}\put(11,11){\line(4,-5){3}}\put(12,12){\line(4,-5){2}}
\put(13,13){\line(4,-5){1}}
\end{picture}
\end{center}
\caption{\scriptsize Strong and weak defining hyperplanes of the PPS; Property 3,}\label{1}
\end{figure}

The PPS of the CCR model is depicted in Figure (\ref{1}). In Figure (1), DMUs $D_1$ and $D_2$ are extreme DMUs and CCR-efficient DMU $D_3$, that lies on the strong defining hyperplane $H_1$ is non-extreme DMUs.

The input-oriented CCR model, corresponds to $DMU_{k}$, $k\in J$, is given by:
\begin{equation} \label{CCRi}
\begin{array}{rlllc}
\min&\displaystyle\ \theta-\epsilon(\sum_{i=1}^{m}s_{i}^{-}+\sum_{r=1}^{s}s_{r}^{+})\\
s.t.&\displaystyle\sum_{j\in J}\lambda_{j}y_{rj}-s_{r}^{+}=y_{rk},&r=1,...,s\\
&\displaystyle\sum_{j\in J}\lambda_{j}x_{ij}+s_{i}^{-}=\theta x_{ik},&i=1,...,m\\
&\lambda_{j}\geq0,&{j\in J}\\
&s_{i}^{-}\geq0,&i=1,...,m\\
&s_{r}^{+}\geq0,&r=1,...,s\\ &\theta&free
\end{array}
\end{equation}Also, the output-oriented CCR model, corresponds to $DMU_{k}$, $k\in
J$, is as follows:\begin{equation} \label{CCRo}
\begin{array}{rlllc}
\max&\displaystyle\ \varphi+\epsilon(\sum_{i=1}^{m}t_{i}^{-}+\sum_{r=1}^{s}t_{r}^{+})\\
s.t.&\displaystyle\sum_{j\in J}\lambda_{j}y_{rj}-t_{r}^{+}=\varphi y_{rk},&r=1,...,s\\
&\displaystyle\sum_{j\in J}\lambda_{j}x_{ij}+t_{i}^{-}= x_{ik},&i=1,...,m\\
&\lambda_{j}\geq0,&{j\in J}\\
&t_{i}^{-}\geq0,&i=1,...,m\\
&t_{r}^{+}\geq0,&r=1,...,s\\
&\varphi&free
\end{array}
\end{equation}
where $\epsilon$ is non-Archimedean small and positive number. Models (\ref{CCRi}) and (\ref{CCRo}) are called envelopment forms (with non-Archimedean number).\\
$DMU_k$ is said to be \textit{strong efficient} (\textit{CCR-efficient}) if and only if either ({\it i}) or ({\it ii}) happen:\
\begin{itemize}
\item[({\it i})] $\theta^*=1$ and $(\textbf{s}^{+*}$, $\textbf{s}^{-*})=(\textbf{0, 0})$
\item[({\it ii})] $\varphi^*=1$ and $(\textbf{t}^{+*}$, $\textbf{t}^{-*})=(\textbf{0, 0})$
\end{itemize}
$DMU_k$ is said to be \textit{weak efficient} if and only if either ({\it v}) or ({\it iv}) happen:\
\begin{itemize}
\item[({\it v})] $\theta^*=1$ and $(\textbf{s}^{+*}$, $\textbf{s}^{-*})\neq (\textbf{0, 0})$
\item[({\it iv})] $\varphi^*=1$ and $(\textbf{t}^{+*}$, $\textbf{t}^{-*})\neq(\textbf{0, 0})$
\end{itemize}
Note that if $\theta^*<1$ and $\varphi^*>1$ then $DMU_k$ is an interior point of the PPS. \footnote{(*)
is used for optimal solution.}%---------------
\maketitle  \textit{{\scriptsize}}
\\Each interior DMU and weak efficient DMU in the CCR model is said to be a \textit{CCR-inefficient} $DMU$.\\{\it Efficient Frontier} is the set of all points (real or virtual DMUs) with efficiency score is equal to unity ($\theta^*=1$ or $\varphi^*=1$).\\Efficient frontier is divided into two categories:

$i)$ \textit{Strong efficient frontier} is the set of all (real or virtual) strong efficient (CCR efficient)

DMU.

$ii)$ \textit{Weak efficient frontier} in which all it's relative interior points (real or virtual DMUs),

are weak efficient DMUs.\\$DMU_k=(X_k, Y_k)$ is said to be \textit{non-dominated} if and only if there is not any $DMU=(X, Y)$ (real or virtual) such that:\\ $(-X_k, Y_k)\geq(-X, Y)$ and $(-X_k, Y_k)\neq(-X, Y)$.\\\\We use the following theorem in the next section.\\
{\bf Theorem 1:}\hskip 0.3cm{\it There does not exist any virtual DMU (a member of the PPS) that dominates an DEA-efficient DMU.}\\
\textbf{Proof.} See H. Fukuyama et. al. (2012).\\\\
The dual of models (\ref{CCRi}) and (\ref{CCRo}) (without $\epsilon$ i.e. $\epsilon$=0), which are called multiplier forms, are as models (\ref{multii}) and (\ref{multio}), respectively:
\begin{equation} \label{multii}
\begin{array}{rlllc}
\max&\displaystyle\sum_{r=1}^{s}u_{r}y_{rk}\\
s.t.&\displaystyle\sum_{r=1}^{s}u_{r}y_{rj}-\sum_{i=1}^{m}v_{i}x_{ij}\leq
0,&j=1,...n\\
&\displaystyle\sum_{i=1}^{m}v_{i}x_{ik}=1,&&\\
&u_{r}\geq0,&r=1,...,s\\
&v_{i}\geq0,&i=1,...,m\\ \end{array}\end{equation}
\begin{equation} \label{multio}
\begin{array}{rlllc}
\min&\displaystyle\sum_{i=1}^{m}v_{i}x_{ik}\\
s.t.&\displaystyle\sum_{i=1}^{m}v_{i}x_{ij}-\sum_{r=1}^{s}u_{r}y_{rj}\geq
0,&j=1,...n\\
&\displaystyle\sum_{r=1}^{s}u_{r}y_{rk}=1,&&\\
&u_{r}\geq0,&r=1,...,s\\
&v_{i}\geq0,&i=1,...,m.\\
\end{array}\end{equation}

$DMU_k$ is strong efficient if there exists at least one optimal solution $(u^*, v^*)$ for (\ref{multii}) with $(u^*, v^*)>\mathbf{0}$, and $u^*y_{k}=1$ in which $u^*=(u_1^*, u_2^*, ..., u_s^*)$ and $v^*=(v_1^*, v_2^*,..., v_m^*)$. Also $DMU_k$ is weak efficient if $u^*y_{k}=1$ and no $(u^*, v^*)>\mathbf{0}$ exists. In this case there exist at least one $r$ (or $i$) so that $u^*_r=0$ (or $v^*_i=0$) in all optimal solution of model (\ref{multii}) or (\ref{multio}) (see example 4.3.). In Figure (\ref{1}), DMUs $D_1$, $D_2$ and $D_3$ are strong efficient and $D'_2$ is weak efficient DMU. The evaluation of $D_3$ and $D'_2$ shows that, model (\ref{multii}) has unique optimal solution, which defines two supporting\footnote{For definition and properties see Bazaraa et al. (1990).} defining hyperplanes $H_1$ and $H_2$ passing through $D_3$ and $D'_2$, respectively. On the other hand, the evaluation of $D_1$ indicate that, model (\ref{multii}) has alternative optimal solutions, which defines an infinite number of supporting hyperplanes passing through $D_1$. Only two of these hyperplanes (i.e. $H_1$ and $H_3$) are defining hyperplanes. In fact, if $(u^*, v^*)$ is an \textit{unique} optimal solution of model (\ref{multii}) then $u^{t*}y-v^{t*}x=0$ is the equation of \textit{defining} hyperplane of the PPS. In addition, if $(u^*, v^*)>0$, $u^{t*}y-v^{t*}x=0$ is the equation of \textit{strong} defining hyperplane of the PPS (see Definition 1). Otherwise, if some components of $(u^*, v^*)$ are zero, then $u^{t*}y-v^{t*}x=0$ is the equation of \textit{weak} defining hyperplane of the PPS (see Definition 2). A similar discussion holds for model (\ref{multio}).\\ Note that the hyperplanes $H_2$ and $H_3$ are weak defining hyperplanes and $H_1$ is strong defining hyperplane of the PPS of the CCR (see Definitions 1 and 2.).

In this paper, corresponding to each strong efficient DMU $DMU_j=(x_{1j}, ..., x_{mj}, y_{1j}, ..., y_{sj})$ we consider virtual DMUs $DMU'_j=(x_{1j}, ..., x_{lj}+\alpha, ..., x_{mj}, y_{1j}, ..., y_{sj})$ and $DMU''_j=(x_{1j}, ..., x_{mj},\\ y_{1j}, ..., y_{qj}-\gamma, ..., y_{sj})$, in which $\alpha, \gamma>0$. These virtual DMUs are either interior point of the PPS of the CCR model or lie on the some weak defining hyperplanes (see Definition 2 and properties 2-7). In the latter case we call these virtual DMUs as ``{\it weak efficient virtual DMU}'', hereafter. (See $DMU$ $D^{'}_{2}$ in Figure \ref{1} and $DMUs$ $D^{'}_{1}$, $D^{'''}_{1}$ in Figure \ref{2}, for example).\\\\
\textbf{Definition 1.} {\it The supporting hyperplane $H=\{(x,y)|\ \bar{u}^{t}y- \bar{v}^{t}x=0, (\bar{u}, \bar{v})\geq0, (\bar{u}, \bar{v})\neq0\}$ of the PPS of the CCR model is strong defining hyperplane of the PPS if only if it is ``defining" and $m+s-1$(=the number of outputs and inputs minus one) strong efficient DMUs of the PPS, which are linear independent, lie on $H$. (In this case, all components of $(\bar{u}, \bar{v})$ are positive.)}\\\\
\textbf{Definition 2.} {\it The supporting hyperplane $H=\{(x,y)|\ \bar{u}^{t}y- \bar{v}^{t}x=0, (\bar{u}, \bar{v})\geq0, (\bar{u}, \bar{v})\neq0\}$ of the PPS of the CCR model is weak defining hyperplane of the PPS if and only if it is ``defining" and $m+s-1$ weak efficient virtual and strong efficient DMUs of the PPS, which are linear independent, lie on H. (In this case, some components of $(\bar{u}, \bar{v})$ are zero.)}\\\\
\textbf{Remark 1}: In the equation of weak defining hyperplane, if $\bar{u}_q=0$ (or $\bar{v}_l$=0), then, this hyperplane is vertical to hyperplane $y_q=0$ (or $x_l=0$). In the case of $x_l=0$, the weak defining hyperplane passes through of $l^{th}$ axis of input.\\\\
\textbf{Remark 2}: If the number of strong efficient DMUs are less than $m+s-1$ then all defining hyperplanes of the PPS are weak defining hyperplanes (because, by Definition 1, at least $m+s-1$ strong efficient DMUs are needed to construct strong defining hyperplane).

In this research, we first find the extreme DMUs of the PPS of the CCR model, lying on the some weak defining hyperplanes, and then using models (\ref{tak}) and (\ref{tao}), the foregone weak efficient virtual DMUs are found. By using them, we find the weak defining hyperplane of the PPS of the CCR model.

{\it Throughout this paper, we must assume that there are not any two strong efficient DMUs as $(x, y)$ and $(tx, ty)$ for all $t>0$ and $t\neq1$. Otherwise, one of them must be deleted.}
\section{Identifying equations of weak defining hyperplanes}\
In this section, we identify the equations of weak defining hyperplanes of the PPS of the CCR model in the following way. First, we evaluate each $DMU_{k},~(k\in J)$ using, models (\ref{CCRi}) or (\ref{CCRo}). Then, we hold all CCR-efficient DMUs, and remove other DMUs. Suppose that the set of all CCR-efficient DMUs is denoted by $E$. Corresponding to each $DMU_{k}=(x_{1k}, ..., x_{mk}, y_{1k}, ..., y_{sk}),~(k\in E)$, we solve the following models:
\begin{equation}
\begin{array}{rlllc}\label{tak}
\min&\displaystyle\ \theta_{l}^{k}\\\\
s.t.&\displaystyle\sum_{j\in E-\{k\}}\lambda_{j}^{k}x_{lj}\leq\theta_{l}^{k} x_{lk}\\
&\displaystyle\sum_{j\in E-\{k\}}\lambda_{j}^{k}x_{ij}\leq x_{ik},&i=1,...,m&i\neq l\\
&\displaystyle\sum_{j\in E-\{k\}}\lambda_{j}^{k}y_{rj}\geq y_{rk},&r=1,...,s&\\
&\lambda_{j}^{k}\geq0,&j\in E-\{k\}\\
&\theta_{l}^{k}&free &l=1,...,m&
\end{array}
\end{equation}
\begin{equation} \label{tao}
\begin{array}{rlllc}
\max&\displaystyle\ \varphi_{q}^{k}\\\\
s.t.&\displaystyle\sum_{j\in E-\{k\}}\mu_{j}^{k}x_{ij}\leq x_{ik},&i=1,...,m&\\
&\displaystyle\sum_{j\in E-\{k\}}\mu_{j}^{k}y_{qj}\geq\varphi_{q}^{k}y_{qk},&&\\
&\displaystyle\sum_{j\in E-\{k\}}\mu_{j}^{k}y_{rj}\geq y_{rk},&r=1,...,s&r\neq q\\
&\mu_{j}^{k}\geq0,&j\in E-\{k\}\\
&\varphi_{q}^{k}& free&q=1,...,s&
\end{array}
\end{equation}
The following properties hold for models (\ref{tak}) and (\ref{tao}). By property 1 we can find all extreme CCR-efficient DMUs. The properties 3, 4, 6, and 7 provide the necessary and sufficient conditions for lying an extreme CCR-efficient DMU on the weak defining hyperplane.\\\\
{\bf Property 1:}\hskip 0.3cm{\it In model (\ref{tak}) (or (\ref{tao})), if for some $l$ (or $q$), $\theta_{l}^{k*}>1$ (or $\varphi_{q}^{k*}<1$) or if for some $l$(or $q$), model (\ref{tak}) (or model (\ref{tao})) is infeasible, then, $DMU_{k}$ is an extreme DMU and vice versa.}\\
\textbf{Proof.} Suppose that $\theta_{l}^{k*}>1$. First, we show that $DMU_k$ is CCR-efficient. By contradiction let $DMU_k$ is inefficient. At optimality of model (\ref{CCRi}), two cases are happened:
\begin{itemize}
\item[{\it (i)}] $\theta^*=1$ and $(\textbf{s}^{+*}, \textbf{s}^{-*})\neq\textbf{0}$
\item[{\it (ii)}] $\theta^*<1$
\end{itemize}
in each cases it can be shown that $\theta_{l}^{k*}\leq1$, a contradiction.\\ Now we show that $DMU_k$ is ,in fact, an extreme CCR-efficient DMU. By contradiction suppose that $DMU_k$ is an non-extreme CCR-efficient. So, the following system has solution:\\
\begin{equation}
\begin{array}{rlllc}\label{tau}
\displaystyle\sum_{j\in E'}\lambda_{j}x_{j}=x_{k},\\
\displaystyle\sum_{j\in E'}\lambda_{j}y_{j}=y_{k},\\
\lambda_{j}\geq0,{j\in E'}
\end{array}
\end{equation}
Suppose that $(\bar{\lambda}_{j}, j\in E')$ is a solution of the above system. If $\bar{\lambda}_{j}=0$ then, $(\theta_{l}^{k}=1, \lambda_{j}=\bar{\lambda}_{j}, j\in E'-\{k\})$ is a solution of model (\ref{tak}). Therefor, $\theta_{l}^{k*}\leq1$, a contradiction. On the other hand if $\bar{\lambda}_{j}\neq0$, we rewrite system (\ref{tau}) as follows:\\
\begin{equation*}
\begin{array}{rlllc}
\displaystyle\sum_{j\in E'-\{k\}}\bar{\lambda}_{j}x_{j}=(1-\bar{\lambda}_{j})x_{k},\\
\displaystyle\sum_{j\in E'-\{k\}}\bar{\lambda}_{j}y_{j}=(1-\bar{\lambda}_{j})y_{k},\\
\end{array}
\end{equation*}
By divided both side of the above equations by $(1-\bar{\lambda}_{j}>0)$; we obtain a solution of model (\ref{tak}) as $(\theta_{l}^{k}=1, \lambda_{j}=\frac{\bar{\lambda}_{j}}{1-\bar{\lambda}_{j}}, j\in E'-\{k\})$. Therefor, $\theta_{l}^{k*}\leq1$, a contradiction. Thus, $DMU_k$ is an extreme CCR-efficient DMU. Now, suppose that for some $l$, model (\ref{tak}) is infeasible. In the similar manner, it can be shown that $DMU_k$ is an extreme CCR-efficient DMU. Conversely, suppose that $DMU_k$ is extreme DMU and model (\ref{tak}) is feasible. We show that  $\theta_{l}^{k*}\geq1$. Consider the following corresponding to $DMU_{k}$:
\begin{equation} \label{takk}
\begin{array}{rlllc}
\min&\displaystyle\ \theta_{l}^{'k}\\
s.t.&\displaystyle\sum_{j\in E'}\lambda_{j}^{k}x_{lj}+s_{l}^{-}=\theta_{l}^{'k} x_{lk}\\
&\displaystyle\sum_{j\in E'}\lambda_{j}^{k}x_{ij}+s_{i}^{-}=x_{ik},&i=1,...,m&i\neq l \\
&\displaystyle\sum_{j\in E'}\lambda_{j}^{k}y_{rj}-s_{r}^{+}=y_{rk},&r=1,...,s&\\
&\lambda_{j}^{k}\geq0,&j\in E'\\
&s_{i}^{-}\geq0, s_{r}^{+}\geq0&i=1,...,m,
&r=1,...,s\\
&\theta_{l}^{'k}&free &l=1,...,m&
\end{array}
\end{equation}
Now suppose that $\theta^{*}(=1)$, $\theta_{l}^{'k*}$ and $\theta_{l}^{k*}$ are the optimal objective functions of the models (\ref{CCRi}), (\ref{takk}) and (\ref{tak}) with respect to $DMU_{k}$, respectively. It is not difficult to show that $\theta^{*}\leq\theta_{l}^{'k*}\leq\theta_{l}^{k*}$. Therefor, $\theta_{l}^{k*}\geq1$. This completes the proof.\hspace{8cm}$\square$\\
{\bf Corollary:}\hskip 0.3cm{\it In models (\ref{tak}) and (\ref{tao}), for each $l$ and $q$, $\theta_{l}^{k*}=\varphi_{q}^{k*}=1$ if and only if $DMU_{k}$ is a non-extreme CCR-efficient DMU.}\\
\textbf{Proof.} Omitted.\\\\
%\textbf{Property 2}: In model (\ref{tak}) (or (\ref{tao})), if for some $l$ (or $q$), $\theta_{l}^{k*}>1$ (or $\varphi_{q}^{k*}<1$) or if for some $l$(or $q$), model (\ref{tak})
%(or model (\ref{tao})) is infeasible, then, strong efficient $DMU_{k}$ lies on the extreme ray (edge) of PPS and vice versa.\\\\
\textbf{Property 2}: In a single input case, for each $DMU_k=(x_{1k}, y_{1k}, ..., y_{sk})$, the virtual DMU $DMU^{'}_k=(x_{1k}+\alpha, y_{1k}, ..., y_{sk})$, in which $\alpha>0$, is an interior point of the PPS of the CCR model.\\\\
\textbf{Proof.} First, we add $DMU^{'}_k$ to the PPS and then, evaluate its performance by the input and output-oriented CCR model \big(see models (\ref{CCRi}) and (\ref{CCRo})\big). It is enough to show that $\theta^*<1$ and $\varphi^*>1$. Consider the input-oriented CCR model corresponding to virtual DMU $DMU^{'}_k$ as follows:
\begin{equation} \label{p3}
\begin{array}{rlllc}
\min&\displaystyle\ \theta\\\\
s.t.
&\displaystyle\sum_{j\in E}\lambda_{j}x_{1j}+\mu_{k}(x_{1k}+\alpha)\leq \theta(x_{1k}+\alpha),&\\
&\displaystyle\sum_{j\in E}\lambda_{j}y_{rj}+\mu_ky_{rk}\geq y_{rk},&r=1,...,s\\
&\lambda_{j}\geq0,&{j\in E}\\
&\theta&free
\end{array}
\end{equation}
$\big(\bar{\lambda}_{j}=0 (j\neq k), \bar{\lambda}_{k}=1, \bar{\mu}_k=0, \bar{\theta}=\frac{x_{1k}}{x_{1k}+\alpha}(<1)\big)$ is a feasible solution of (\ref{p3}). Since model (\ref{p3}) has a minimization-type objective function, $\theta^*<1$; where ``*" is used to indicate optimality. In a similar manner, it can be shown that in output-oriented maximization problem, $\varphi^*>1$. Therefore, $DMU_k$ is an interior point of PPS. This completes the proof.\hspace{6cm}$\square$\\

\begin{figure}
[b]
\begin{center}
\unitlength=2.5mm
\begin{picture}(18,2)(-4,-4)
\thicklines\put(0,0){\vector(1,0){14}}
\put(0,0){\vector(0,1){14}}
\thicklines\put(0,0){\vector(-1,-1){7}}
\put(0,0){\makebox(31,0)[]{\shortstack[b]{$y_{2}$}}}
\put(0,0){\makebox(0,30)[]{\shortstack[b]{x}}}
\put(0,0){\makebox(-15,-15)[]{\shortstack[b]{$y_{1}$}}}
\put(2,6){{\circle*{.6}}\makebox(2,-2)[]{\shortstack[b]{$D_{1}$}}}
\put(3,7){{\circle*{.6}}\makebox(2,-1)[]{\shortstack[b]{$D'''_{1}$}}}
\put(0.5,6){{\circle*{.6}}\makebox(-6,-2)[]{\shortstack[b]{$D'_{1}$}}}
\put(2,9.5){{\circle{.4}}\makebox(0,1.5)[]{\shortstack[b]{\small$D''_{1}$}}}
\put(0,0){{\circle*{.6}}\makebox(0,-2)[]{\shortstack[b]{O}}}
\thicklines\put(0,0){\line(1,3){6}}
\thinlines\put(0,0){\line(-1,5){3}}
\thinlines\put(0,0){\line(1,2){9}}
\thinlines\put(0.3,1){\line(-1,0){.5}}\put(0.6,2){\line(-1,0){1}}\put(0.9,3){\line(-1,0){1.5}}\put(1.3,4){\line(-1,0){2}}\put(1.7,5){\line(-1,0){2.5}}\put(2.,6){\line(-1,0){3}}
\put(2.3,7){\line(-1,0){3.7}}\put(2.7,8){\line(-1,0){4.3}}
\thinlines\put(3,9){\line(-1,0){4.8}}\put(3.7,11){\line(-1,0){5.9}}\put(4,12){\line(-1,0){6.4}}
\put(0.3,1){\line(1,1){1.1}}\put(0.6,2){\line(1,1){1.18}}\put(1,3){\line(1,1){1.18}}\put(1.3,4){\line(1,1){1.3}}\put(1.6,5){\line(1,1){1.8}}
\put(2,6){\line(1,1){2}}\put(2.4,7){\line(1,1){2.3}}\put(2.7,8){\line(1,1){2.5}}\put(3,9){\line(1,1){3}}\put(3.3,10){\line(1,1){3.4}}
\put(3.6,11){\line(1,1){3.8}}\put(4,12){\line(1,1){4.1}}
\end{picture}
\end{center}
\caption{\scriptsize Properties 2 and 6.}\label{2}
\end{figure}
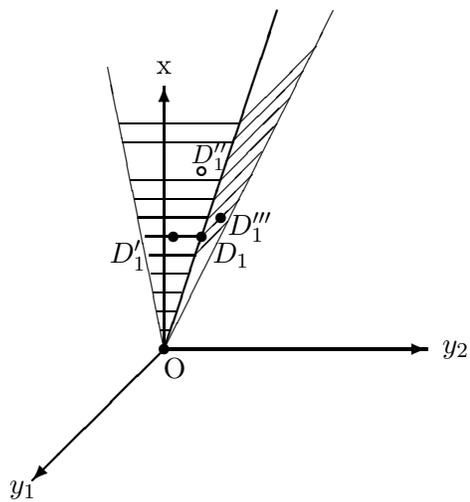

In Figure \ref{2}, corresponding to DMU $D_1=(x_{11}, y_{11}, y_{21})$, virtual DMU $D''_1=(x_{11}+\alpha, y_{11}, y_{21})$ is an interior point of the PPS.\\\\
\textbf{Property 3}: In a multiple inputs case, if for some $l$, model (\ref{tak}) is infeasible, then CCR-efficient $DMU_k$ lies on the weak defining hyperplane, which passes through the $l$th axis of input.\\\\
\textbf{Proof.} We show that if for one $l$, model (\ref{tak}) is infeasible, then, virtual DMU $DMU^{'}_k=(x_{1k}, ..., x_{(l-1)k}, x_{lk}+\alpha, x_{(l+1)k}, ..., x_{mk}, y_{1k}, ..., y_{sk})$, in which $\alpha>0$, is on the weak defining hyperplane, which passes through the $l$th axis of input. For this aim, we show that in the performance evaluation of $DMU^{'}_k$ using model (\ref{CCRo}); $\varphi^*=1$. Consider model (\ref{CCRo}) corresponding to virtual DMU $DMU^{'}_k$ as follows (without $\epsilon$):
\begin{equation} \label{10}
\begin{array}{rlllc}
\max&\displaystyle\ \varphi\\\\
s.t.&\displaystyle\sum_{j\in E}\lambda_{j}y_{rj}+\mu_{k}y_{rk}\geq\varphi y_{rk},&r=1,...,s\\
&\displaystyle\sum_{j\in E}\lambda_{j}x_{ij}+\mu_{k}x_{ik}\leq x_{ik},&i=1,...,m,&i\neq l\\
&\displaystyle\sum_{j\in E}\lambda_{j}x_{lj}+\mu_{k}(x_{lk}+\alpha)\leq x_{lk}+\alpha\\
&\mu_k, \lambda_{j}\geq0,&{j\in E}\\
&\varphi&free
\end{array}
\end{equation}
By contradiction, suppose that $\big(\lambda^{*}_{j}\  (j\in E),\ \mu^{*}_{k},\ \varphi^*(>1)\big)$ is the optimal solution of (\ref{10}). The constraints of model (\ref{10}) can be written as follows:
\begin{equation} \label{11}
\begin{array}{rlllc}
&\displaystyle\sum_{j\in E-\{k\}}\lambda_{j}^{*}y_{rj}>(1-\lambda_{k}^{*}-\mu_{k}^{*})y_{rk},&r=1,...,s\\
&\displaystyle\sum_{j\in E-\{k\}}\lambda_{j}^{*}x_{ij}\leq(1-\lambda_{k}^{*}-\mu_{k}^{*}) x_{ik},&i=1,...,m,&i\neq l\\
&\displaystyle\sum_{j\in E-\{k\}}\lambda_{j}^{*}x_{lj}\leq (1-\lambda_{k}^{*}-\mu_{k}^{*})x_{lk}+(1-\mu_{k}^{*})\alpha\\
\end{array}
\end{equation}
From model (\ref{11}), it is easy to show that $1-\lambda_{k}^{*}-\mu_{k}^{*}>0$. Divide both sides of model (\ref{11}) by $1-\lambda_{k}^{*}-\mu_{k}^{*}>0$ and define $\bar{\mu}_{j}=\displaystyle\frac{\lambda_{j}^{*}}{1-\lambda_{k}^{*}-\mu_{k}^{*}}$, $j\in E-\{k\}$; so, model (\ref{11}) becomes as follows:
\begin{equation} \label{12}
\begin{array}{rlllc}
&\displaystyle\sum_{j\in E-\{k\}}\bar{\mu}_{j}y_{rj}>y_{rk},&r=1,...,s\\
&\displaystyle\sum_{j\in E-\{k\}}\bar{\mu}_{j}x_{ij}\leq x_{ik},&i=1,...,m,&i\neq l\\
&\displaystyle\sum_{j\in E-\{k\}}\bar{\mu}_{j}x_{lj}\leq x_{lk}+\beta\\
\end{array}
\end{equation}
in which $\beta=\Big(\displaystyle\frac{1-\mu_{k}^{*}}{1-\lambda_{k}^{*}-\mu_{k}^{*}}\Big)\alpha$. Since $\beta>0$, there is $\hat{\theta}>0$ so that $x_{lk}+\beta=\hat{\theta}x_{lk}$; therefore, the constraints of model (\ref{11}) can be rewritten as follows:
\begin{equation*}
\begin{array}{rlllc}
&\displaystyle\sum_{j\in E-\{k\}}\bar{\mu}_{j}y_{rj}>y_{rk},&r=1,...,s\\
&\displaystyle\sum_{j\in E-\{k\}}\bar{\mu}_{j}x_{ij}\leq x_{ik},&i=1,...,m,&i\neq l\\
&\displaystyle\sum_{j\in E-\{k\}}\bar{\mu}_{j}x_{lj}\leq \hat{\theta}x_{lk}\\
\end{array}
\end{equation*}
So, $(\bar{\mu_j}\ \ (j\in E-\{k\}), \hat{\theta})$ is a feasible solution for model (\ref{tak}); a contradiction. This implies that $\varphi^*=1$ i.e. $DMU'_k$ lies on the efficient frontier. Now, since $DMU'_k$ is dominated by CCR-efficient $DMU_k$, so, $DMU'_k$ lies on the weak efficient frontier (hyperplane). Moreover, it is easy to shows that in the equation of this weak efficient hyperplane; $v_l=0$ and so by remark 1 this hyperplane passes through the $l$th axis of input. This completes the proof.\hspace{2.5cm}$\square$\\

In Figure \ref{1}, model (\ref{tak}) corresponding to DMU $D_2=(x_{12}, x_{22}, y_{2})$ with $l=1$ is infeasible; so, virtual DMU $D'_2=(x_{12}+\alpha, x_{22}, y_{2})$ is on the weak defining hyperplane, which passes through $x_1$-axis and vertical to hyperplane $x_1$=0.\\\\
The following property is, in fact, the converse of property 3.\\
\textbf{Property 4}: In a multiple inputs case, if extreme CCR-efficiency DMU $DMU_k=(x_{1k}, ..., x_{lk}, ...\\, x_{mk}, y_{1k}, ..., y_{sk})$ lies on the weak defining hyperplane which passes through the $l$th axis of input (vertical to hyperplane $x_l$=0); then model (\ref{tak}) is infeasible.\\\\
\textbf{Proof.} By contradiction, suppose that the model (\ref{tak}) is feasible. The first constraint of the model (\ref{tak}) implies that the optimal solution of the model (\ref{tak}) is bounded. Suppose that, $\big(\theta^{k*}_l, \lambda^*_j~~ (j\neq k)\big)$ is an optimal solution of it. Note that the first constraint of the model (\ref{tak}) is tight at optimality. We first show that $\theta^{k*}_l>1$. By contradiction, suppose that $\theta^{k*}_l\leq1$. If $\theta^{k*}_l<1$ we have:
\begin{equation}
\begin{array}{rlllc}\label{tk}
&\displaystyle\sum_{j\in E-\{k\}}\lambda_{j}^{k*}x_{lj}=\theta_{l}^{k*} x_{lk}<x_{lk}\\
&\displaystyle\sum_{j\in E-\{k\}}\lambda_{j}^{k*}x_{ij}\leq x_{ik},&i=1,...,m&i\neq l\\
&\displaystyle\sum_{j\in E-\{k\}}\lambda_{j}^{k*}y_{rj}\geq y_{rk},&r=1,...,s&\\
\end{array}
\end{equation}
It shows that virtual DMU $$(\displaystyle\sum_{j\in E-\{k\}}\lambda_{j}^{k*}x_{lj}, ..., \displaystyle\sum_{j\in E-\{k\}}\lambda_{j}^{k*}x_{lj}, ..., \displaystyle\sum_{j\in E-\{k\}}\lambda_{j}^{k*}x_{mj}, \displaystyle\sum_{j\in E-\{k\}}\lambda_{j}^{k*}y_{1j}, ..., \sum_{j\in E-\{k\}}\lambda_{j}^{k*}y_{sj})$$ dominates the CCR-efficient $DMU_k$, a contradiction (see Theorem 1). Now, if $\theta^{k*}_l=1$, we have:
\begin{equation}
\begin{array}{rlllc}\label{to}
&\displaystyle\sum_{j\in E-\{k\}}\lambda_{j}^{k*}x_{lj}=x_{lk}\\
&\displaystyle\sum_{j\in E-\{k\}}\lambda_{j}^{k*}x_{ij}\leq x_{ik},&i=1,...,m&i\neq l\\
&\displaystyle\sum_{j\in E-\{k\}}\lambda_{j}^{k*}y_{rj}\geq y_{rk},&r=1,...,s&\\
\end{array}
\end{equation}
At least one of the inequality constraints of (\ref{to}) is a strict inequality, because, otherwise, the CCR-efficient $DMU_k$, is not extreme DMU. So, $\theta^{k*}_l>1$. Therefor, there exist $\beta>0$ so that, $\theta^{k*}_lx_{lk}=x_{lk}+\beta$. This means that, the virtual DMU $DMU^{'}_k=(x_{1k}, ..., x_{(l-1)k}, x_{lk}+\beta, x_{(l+1)k}, ..., x_{mk}, y_{1k}, ..., y_{sk})$ is, in fact, an observed DMU belongs to the PPS of the CCR model. This is a contradiction. Because, we had been eliminated all the CCR-inefficient DMUs from the PPS of the CCR model. The proof is completed.\hspace{10.5cm}$\square$\\\\
\textbf{Property 5}: In a single output case, for each $DMU_k=(x_{1k}, ..., x_{mk}, y_{1k})$, virtual DMU $DMU^{'}_k=(x_{1k}, ..., x_{mk}, y_{1k}-\gamma)$, in which $\gamma>0$, is an interior point of the PPS of the CCR model.\\\\
\textbf{Proof.} The proof is similar to property 2 and so, we omit the details. \hspace{4.3cm}$\square$\\

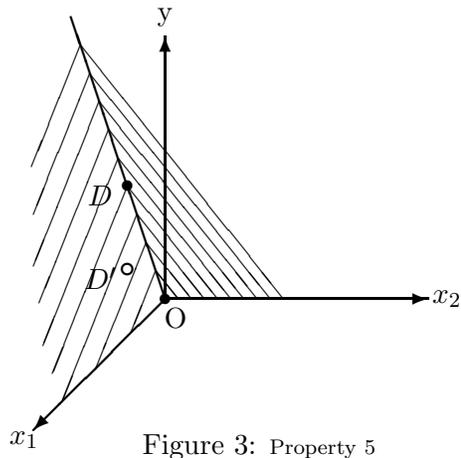
\begin{figure}
[p]
\begin{center}
\unitlength=2.5mm
\begin{picture}(18,-15)(-4,-4)
\thicklines\put(0,0){\vector(1,0){14}}
\put(0,0){\vector(0,1){14}}
\thicklines\put(0,0){\vector(-1,-1){7}}
\put(0,0){\makebox(30,0)[]{\shortstack[b]{$x_{2}$}}}
\put(0,0){\makebox(0,30)[]{\shortstack[b]{y}}}
\put(0,0){\makebox(-15,-15)[]{\shortstack[b]{$x_{1}$}}}
\put(-2,6){{\circle*{.6}}\makebox(-4,-1)[]{\shortstack[b]{$D$}}}
%\put(-1,13){\makebox(-4,-1)[]{\shortstack[b]{$H_1$}}}
\put(-2,1.6){{\circle{.6}}\makebox(-4,-1)[]{\shortstack[b]{$D'$}}}
\put(0,0){{\circle*{.6}}\makebox(0,-2)[]{\shortstack[b]{O}}}
%\thicklines\put(0,0){\line(1,1){15}}
\thicklines\put(0,0){\line(-1,3){5}}
\thinlines\put(-.5,1.5){\line(-2,-5){1.3}}\put(-1,3){\line(-2,-5){2.7}}\put(-1.5,4.5){\line(-2,-5){4}}\put(-2,6){\line(-2,-5){4.}}
\thinlines\put(-2.5,7.5){\line(-2,-5){4}}\put(-3,9){\line(-2,-5){3.7}}\put(-3.5,10.5){\line(-2,-5){3.5}}\put(-4,12){\line(-2,-5){3.}}\put(-4.5,13.5){\line(-2,-5){2.6}}
\thinlines\put(-.5,1.5){\line(4,-5){1.2}}\put(-1,3){\line(4,-5){2.4}}\put(-1.5,4.5){\line(4,-5){3.6}}\put(-2,6){\line(4,-5){4.8}}\put(-2.5,7.5){\line(4,-5){6}}
\thinlines\put(-3,9){\line(4,-5){7.2}}\put(-3.5,10.5){\line(4,-5){8.4}}\put(-4,12){\line(4,-5){9.6}}\put(-4.5,13.5){\line(4,-5){10.8}}
\end{picture}
\caption{\scriptsize Property 5}\label{3}
\end{center}
\end{figure}

In Figure \ref{3}, virtual DMU $D'=(x_{1}, x_{2}, y_{1}-\gamma)$, corresponding to DMU $D=(x_{1}, x_{2}, y_{1})$ is an interior point of PPS).\\\\
\textbf{Property 6}: In a multiple outputs case, if for some $q$, model (\ref{tao}) is infeasible, then, CCR-efficient $DMU_k$ lies on the weak defining hyperplane of the PPS, vertical to hyperplane $y_q$=0.\\\\
\textbf{Proof.} The proof is similar to property 3 except it can be shown that in the performance evaluation of $DMU^{'}_k$ using model (\ref{CCRi}); $\theta^*=1$.\hspace{9cm}$\square$\\

In Figure \ref{2}, model (\ref{tao}) corresponding to DMU $D_1=(x_{11}, y_{11}, y_{21})$, with $q=2$, is infeasible, so, virtual DMU $D'_1=(x_{11}, y_{11}, y_{21}-\gamma)$ lies on the weak defining hyperplane of the PPS vertical to hyperplane $y_2$=0.\\\\
The following property is, in fact, the converse of property 6.\\\\
\textbf{Property 7}: In a multiple outputs case, let extreme CCR-efficiency DMU $DMU_k=(x_{1k}, ..., x_{mk},\\y_{1k}, ..., y_{qk}, ...., y_{sk})$ lies on the weak defining hyperplane vertical to hyperplane $y_q=0$; then model (\ref{tao}) is infeasible.\\\\
\textbf{Proof.} The proof is similar to property 4. So, we omit it.\\

Now,  by property 1 we can find all extreme DMUs of the PPS of the CCR model and by Properties 2, 3 and 4 we can find all weak efficient virtual DMUs as $DMU^{'}_k=(x_{1k}, ..., x_{(l-1)k}, x_{lk}+\alpha, x_{(l+1)k}, ..., x_{mk}, y_{1k}, ..., y_{sk})$, in which $\alpha>0$ and by Properties 5, 6 and 7 we can find all weak efficient virtual DMUs as $DMU'_k=(x_{1k}, ..., x_{mk},  y_{1k}, ..., y_{qk}-\beta, ...., y_{sk})$, in which $\beta>0$. Put indices of the weak efficient virtual DMUs, in $F$. Without lose of generality we can assume that $F\cup E=\{DMU_1, . . . ,DMU_L\}$. Consider the set $G=\{1, ..., L\}$. Now, we use the method proposed by Jahanshahloo et. al (2007) to find all weak defining hyperplanes of the PPS of the CCR model with some modifications. We modify their method for simplifying the process of finding coplanar DMUs. More importantly, we need the following theorems:\\\\
{\bf Theorem 2:}\hskip 0.5cm{\it Let $(x_{p},y_{p})$ and $(x_{q},y_{q})$ be observed DMUs that lie on a strong (or weak) supporting hyperplane, then each convex combination of them is on the same hyperplane.}\\ \textbf{Proof.} See Jahanshahloo et. al (2007).\\\\
{\bf Theorem 3:}\hskip 0.5cm{\it Consider $(x_{p},y_{p})$ and $(x_{q},y_{q})$ are two observed DMUs that lie on different hyperplanes (excluding their intersection, if it is not empty). Then every point (virtual DMU) which is obtained by strict convex combination of them is an interior point of PPS. In other words, this virtual DMU is radial inefficient.}\\
\textbf{Proof.} See Jahanshahloo et. al (2007).\\

The method is as follows:\\
Using the above Theorems we first determine all coplanar DMUs in $G$ (i.e. all DMUs that are on the same defining hyperplane). Take a distinct pair $DMU_{p}$ and $DMU_{q}$, where $p$ and $q$ belong to $G$, and construct a virtual DMU as follows: $DMU_{k}$=$\frac{1}{2}DMU_{p}$+$\frac{1}{2}DMU_{q}$. Using the DEA models, we can determine whether or not $DMU_{k}$ is efficient. In the first case, by Theorem 3, $DMU_{p}$ and $DMU_{q}$ are on the same hyperplane; in the second case, they are not (by Theorem 2). For each $l\in G$, define $G_{l}=\big\{j|\ DMU_{l}\ $and$\ DMU_{j},\ j\in G,\ $are coplanar$\big\}$. It is obvious that if $\{l_1, l_2, ..., l_p\}\subseteq G_{l_1}\cap G_{l_2}\cap...\cap G_{l_p}$, then, DMUs $l_1, l_2, ..., l_p$ are coplanar.\\\\ \textbf{Example.} In Figure (\ref{4}), 1, 2, ..., 6 are six DMUs. Since there exists a plane that binding from $1$ and $2$, therefore $1\in G_{2}$ and $2\in G_{1}$ and so on. Then, $G_{1}=\{1,2,3\}$, $G_{2}=\{1,2,3,4,5\}$, $G_{3}=\{1,2,3,4,5\}$, $G_{4}=\{2,3,4,5,6\}$, $G_{5}=\{2,3,4,5,6\}$, $G_{6}=\{4,5,6\}$.\\ $\{1, 2, 3\}\subseteq G_{1}\cap G_{2}\cap G_{3}$ therefor, DMUs 1, 2, 3, are coplanar.\\ $\{2, 3, 4, 5\}\subseteq G_{2}\cap G_{3}\cap G_4\cap G_{5}$ therefor, DMUs 2, 3, 4, 5 are coplanar.\\ $\{4, 5, 6\}\subseteq G_{4}\cap G_{5}\cap G_{6}$ therefor, DMUs 4, 5, 6 are coplanar.\\

We choose an arbitrary $m+s-1$ members of $G$ such that $\{l_1, l_2, ..., l_{m+s-1}\}\subseteq G_{l_1}\cap G_{l_2}\cap...\cap G_{l_{m+s-1}}$. We call this set $D=\{j_{1},\ldots,j_{m+s-1}\}$.\\\\
Using $D$, a hyperplane can be constructed as follows:\\\\
\[{\scriptsize
\left|\begin{array}{llllll}x_{1}&\cdots&x_{m}&y_{1}&\cdots&y_{s}\\
x_{1j_{1}}&\cdots&x_{mj_{1}}&y_{1j_{1}}&\cdots&y_{sj_{1}}\\.&.&.&.&.&.\\.&.&.&.&.&.\\
.&.&.&.&.&.\\x_{1j_{m+s-1}}&\cdots&x_{mj_{m+s-1}}&y_{1j_{m+s-1}}&\cdots&y_{sj_{m+s-1}}\end{array}
\right|}=0,\hspace{4cm}(12)\]\\\\ where $x_{1},\ldots,x_{m},y_{1},\ldots,y_{s}$ are variables, $x_{pj_{t}}, (p=1,\ldots,m,\ t=1,\ldots,$ $ m+s-1)$ is $p$ $th$ input of $DMU_{j_{t}}$ and $y_{qj_{t}}\ (q=1,\ldots,s;\ t=1,\ldots,m+s-1)$ is $q$ $th$ output of $DMU_{j_{t}}$.\\Suppose that the equation of the above mentioned hyperplane is in the form of $P^{t}z=0$, where $z=(x_{1},\ldots,x_{m},y_{1},\ldots,y_{s})$, and $P$ is the gradient vector of the hyperplane. Considering Theorem 4, we can find all equations of weak and strong defining hyperplanes of PPS.\\\\
{\bf Theorem 4:}\hskip 0.5cm{\it Consider $H=\{z|\ P^{t}z=0\}$ so that $P^{t}z=0$ is constructed by (12). Suppose $w=(x_{1}^{w},\ldots,x_{m}^{w},y_{1}^{w},\ldots,y_{s}^{w})$ is defined as follows:\\ $x^{w}_{i}=\max \{x_{ij}|\ j=1,\ldots,n\}\ i=1,\ldots,m$\\ $y^{w}_{r}=\min\{y_{rj}|\ j=1,\ldots,n\}\ r=1,\ldots,s$\\Call $w$ Negative Ideal (NI) if\\ $P^{t}z_j=0\ j\in D$\\ $P^{t}z_j\leq0\ j\in G-D$\\
$P^{t}w<0$\\then H is supporting.}\\
\textbf{Proof.} See Jahanshahloo et. al (2007).

Now we are in the position to put all together the ingredients of the method.\\\\
\textbf{Summary of finding all Weak Defining Hyperplanes' algorithm}
\begin{itemize}
\item\textbf{Step 1}. Evaluate \textbf{n} DMUs with a suitable form of models (\ref{CCRi}) and, (\ref{CCRo}). Hold all CCR-efficient DMUs and remove other DMUs. Put indices of this strong efficient DMUs in $E$.
\item\textbf{Step 2}. Evaluate each CCR-efficient DMUs with models (\ref{tak}) and (\ref{tao}). (Note that in the single output case we use model (\ref{tak}) and in the single input case we use model (\ref{tao})). Denote the index set of weak efficient virtual DMUs by $F$. Suppose that $F\cup E=\{DMU_1, . . . ,DMU_L\}$ and $G=\{1, ..., L\}$.
\item\textbf{Step 3}. For each $p,q\in G$ that $p\neq q$, evaluate $DMU_{k}=\frac{1}{2}DMU_{p}$+$\frac{1}{2}DMU_{q}$ if it is efficient $p\in G_{q}$ and $q\in G_{p}$.
\item\textbf{Step 4}. For each $j\ (j=1,\ldots,L)$, compute $G_{j}$.
\item\textbf{Step 5}. Choose arbitrary $m+s-1$ members of $G$ such that $\{l_1, l_2, ..., l_{m+s-1}\}\subseteq G_{l_1}\cap G_{l_2}\cap...\cap G_{l_{m+s-1}}$. Call this set as $D=\{j_{1},\ldots,j_{m+s-1}\}$.\\ Construct a hyperplane using equation (12). Suppose that the equation of hyperplane is in the form of $P^{t}z=0$ where $z=(x_{1},\ldots,x_{m},\indent\indent y_{1},\ldots, y_{s})$.
\item\textbf{Step 6}. If $P$ has any component less than or equal to zero go to step 8, else let $w=(x_{1}^{w},\ldots, x_{m}^{w}, y_{1}^{w},\ldots, y_{s}^{w})$ is defined as follows:\\\\
$x^{w}_{i}=\max \{x_{ij}|\ j=1,\ldots,n\},~~~~~ i=1,\ldots,m$\\
$y^{w}_{r}=\min \{y_{rj}|\ j=1,\ldots,n\},~~~~~ r=1,\ldots,s$\\If\\ $P^{t}z_{j}=0,\ \ j\in D$\\
$P^{t}z_{j}\leq0,\ \ j\in G-D$\\
$P^{t}w<0,\ \ $\\then $P^{t}z=0$ is supporting. Otherwise, go to step 8.
\item\textbf{Step 7}. If, at least, one of the $m+s-1$ members of $D$ is a weak efficient virtual DMU, then $P^{t}z=0$ is weak defining hyperplane. Otherwise, it is strong defining hyperplane.
\item\textbf{Step 8}. If another subset of $G$ with $m+s-1$ members can be found, go to step 5, else stop.\end{itemize}
\section{Numerical Examples}
\textit{4.1. Single output case}\\

Table 1 shows data for 4 DMUs with two inputs and one output. Using the CCR model (\ref{CCRi}), the CCR-efficient DMUs are determined to be $D_1$, $D_2$, $D_3$. Remove CCR-inefficient DMU $D_4$ from PPS and solve model (\ref{tak}) corresponding to CCR-efficient DMUs $D_1$, $D_2$ and $D_3$.\\
The following results are yielded:\\By property 1, DMUs $D_1$, $D_2$ and $D_3$ lie on the extreme rays of the PPS. Model (\ref{tak}) corresponding to DMU $D_1$ with $l=2$ and DMU $D_3$ with $l=1$ is infeasible. So, by property 3 weak efficient virtual DMUs $D_5=D'_1=(1, 4+\beta, 3)$ and $D_6=D'_3=(5+\alpha, 1, 5)$ lie on the weak defining hyperplanes which pass thought the $\textit{2}$th and $\textit{1}$th axis of inputs, respectively. For convenience, we let $\alpha=\beta=1$. Therefore $E\cup F=\{D_1, D_2, D_3, D_5, D_6\}$ and $G=\{1, 2, 3, 5, 6\}$ . Note that model (\ref{tak}) corresponding to DMU $D_2$ is feasible. So, by property 4, DMU $D_2$ does not lie on any weak defining hyperplane. Also\\
$G_1=\{1, 2, 5\}$,~~~~$G_2=\{2, 1, 3\}$,~~~~ $G_3=\{3, 2, 6\}$,~~~~ $G_5=\{1, 5\}$,~~~~ $G_6=\{6, 3\}$.\\
$\{1, 5\}\subseteq G_1\cap G_5$,~~~~ $\{1, 2\}\subseteq G_1\cap G_2$,\\$\{2, 3\}\subseteq G_2\cap G_3$,~~~ $\{3, 6\}\subseteq G_3\cap G_6$.

$H_{1}$, the first weak hyperplane, is constructed on $D=\{1, 5\}$.\\

$\left|\begin{array}{lll}x_{1}&x_{2}&y\\
1&4&3\\ 1&5&3\end{array} \right|=0$, that yields $y=3x_1$.\\\\ Note that the conditions of Theorem 4 are held and $H_{1}$ is a weak defining hyperplane including $x_2$-axis and vertical to hyperplane $x_2=0$ (because the weak efficient virtual DMU, $D'_{1}$, lies on $H_{1}$). Here, $w=(6,4,2)$.

$H_{2}$, the second weak hyperplane is constructed on $D=\{3, 6\}$.\\

$\left|\begin{array}{lll}x_{1}&x_{2}&y\\
5&1&5\\ 6&1&5\end{array} \right|=0$, that yields $y=5x_2$.\\\\Similarly, the conditions of Theorem 4 are held, and hence $H_{2}$ is a weak defining hyperplane including $x_1$-axis and vertical to hyperplane $x_1=0$ (because the weak efficient virtual DMU, $D'_{3}$, lies on $H_{2}$).\\ Using the proposed method, it was found that DMUs $D_1$, $D_2$ and $D_3$ rest on edges of the PPS. Also, DMUs $D_1$ and $D_3$ lie on the edge intersection of strong and weak defining hyperplanes of the PPS. Moreover, using property 5 there is no any weak defining hyperplane vertical to hyperplane $y=0$.
\begin{table}
\caption{{\scriptsize Data of Numerical Example 1.}}
{\scriptsize
\begin{tabular}{lllllccccccccccccccccccccccccccrr}
  % after \\ : \hline or \cline{col1-col2} \cline{col3-col4} ...
 %&&&& &&&& &&&& &&&& \\
&&&& &&&& &&&& &&&&\\
\hline
  DMU &&&&&&&D1& &&&&&&&D2&&&&&&&&D3&&&&&&&&D4  \\
\hline
  $x_1$ &&&&&&&1& &&&&&&&2&&&&&&&&5&&&&&&&&6 \\

  $x_2$ &&&&&&&4& &&&&&&&2&&&&&&&&1&&&&&&&&1 \\

  $y$ &&&&&&&3& &&&&&&&4&& &&&&&&5&&&&&&&&2 \\
\hline
\end{tabular}}
\end{table}
\begin{table}
\caption{{\scriptsize Example 2. Multiple input and output.}}
\begin{center}
{\scriptsize
\begin{tabular}{lllllllllllllllllllllllllllllll}
\hline
DMU&&&&&&$D_1$&&&&&&$D_2$&&&&&&$D_3$&&&&&&$D_4$&&&&&&$D_5$\\
\hline
$x_1$&&&&&&2&&&&&&1&&&&&&2&&&&&&4&&&&&&3\\
$x_2$&&&&&&3&&&&&&2&&&&&&2&&&&&&2&&&&&&5\\
$y_1$&&&&&&7&&&&&&3&&&&&&4&&&&&&6&&&&&&5\\
$y_2$&&&&&&4&&&&&&5&&&&&&3&&&&&&1&&&&&&2\\
\hline
\end{tabular}}
\end{center}
\end{table}\\\\
\textit{4.2. Multiple outputs and inputs case}\\

Table 2 shows data for 5 DMUs with two inputs and two outputs. Running model (\ref{CCRi}) (or (\ref{CCRo})) shows that DMUs $D_1$, $D_2$, $D_{4}$ are CCR-efficient and other DMUs are CCR-inefficient. Applying models (\ref{tak}) and (\ref{tao}) to each CCR-efficient DMU produces the results reported in Table 3. In Table 3, ''INFES'' and "FES" denotes "infeasible" and "feasible" respectively. For instance, "INFES" in the first row and the second column means that model (\ref{tak}), corresponding to DMU $D_1$ with $l=2$, is infeasible. So, by property 3, $D^{2}_{1}=(2, 3+\alpha, 7, 4)$ is a weak efficient virtual DMU that lies on a weak defining hyperplane passing through $x_2$-axis. using property 1, it was found that all CCR-efficient DMUs lie on the extreme ray. Using properties 3 and 6 and the information of Table 3, all weak efficient virtual DMUs can be determined (see Table 4). For simplicity, we choose $\alpha=1$. For simplicity, we rename CCR-efficient and weak efficient virtual DMUs as follows:\\\\ $U_1=D_1$, $U_2=D_2$, $U_3=D_4$, $U_4=D_{1}^{2}$, $U_5=D_{1}^{4}$, $U_6=D_{2}^{1}$, $U_7=D_{2}^{2}$, $U_8=D_{2}^{3}$, $U_9=D_{4}^{1}$, $U_{10}=D_{4}^{4}$.\\\\Therefore, we have:\\\\ $E=\{U_1, U_2, U_3\}$,\hspace{2cm} $F=\{U_4, U_5, U_6, U_7, U_8, U_9, U_{10}\}$,\\$G=\{U_1, U_2, U_3, U_4, U_5, U_6, U_7, U_8, U_9, U_{10}\}$,\\
$G_1=\{U_1, U_2, U_3, U_4, U_5, U_7, U_{10}\}$,~~~
$G_2=\{U_1, U_2, U_3, U_4, U_6, U_7, U_8, U_9\}$,\\
$G_3=\{U_1, U_2, U_3, U_5, U_6, U_9, U_{10}\}$,~~~
$G_4=\{U_1, U_2, U_4, U_5, U_7\}$,\\
$G_5=\{U_1, U_3, U_4, U_5, U_{10}\}$,~~~~~~~~~~~~~
$G_6=\{U_2, U_3, U_6, U_8, U_9\}$,\\
$G_7=\{U_1, U_2, U_4, U_7, U_8\}$,~~~~~~~~~~~~~
$G_8=\{U_2, U_6, U_7, U_8\}$,\\
$G_9=\{U_2, U_3, U_6, U_9, U_{10}\}$,~~~~~~~~~~~~~
$G_{10}=\{U_1, U_3, U_5, U_9, U_{10}\}$.

$H_{1}$, the first weak hyperplane, is constructed on $D1=\{1, 2, 4\}$, $D2=\{1, 2, 7\}$, $D3=\{2, 1, 7\}$, $D4=\{2, 4, 7\}$,\\

\small{
$\left|\begin{array}{lllll}x_{1}&x_{2}&y_1&y_2\\
2&3&7&4&\\ 1&2&3&5&\\2&4&7&4&\end{array} \right|$=$\left|\begin{array}{lllll}x_{1}&x_{2}&y_1&y_2\\
2&3&7&4&\\ 2&4&7&4&\\1&3&3&5&\end{array} \right|$=$\left|\begin{array}{lllll}x_{1}&x_{2}&y_1&y_2\\
1&2&3&5&\\ 2&3&7&4&\\1&3&3&5&\end{array} \right|$=$\left|\begin{array}{lllll}x_{1}&x_{2}&y_1&y_2\\
1&2&3&5&\\ 2&4&7&4&\\1&3&3&5&\end{array} \right|
=0$}\\\\\\ that yields $-23x_1+6y_1+y_2=0$.\\

$H_{2}$, the second weak hyperplane, is constructed on $D5=\{1, 10, 3\}$, $D6=\{1, 3, 5\}$, $D7=\{1, 5, 10\}$, $D8=\{3, 5, 10\}$\\

$\left|\begin{array}{lllll}x_{1}&x_{2}&y_1&y_2\\
2&3&7&4&\\ 4&2&6&.5&\\4&2&6&1&\end{array} \right|$=$\left|\begin{array}{lllll}x_{1}&x_{2}&y_1&y_2\\
2&3&7&4&\\ 4&2&6&1&\\2&3&7&3&\end{array} \right|$=$\left|\begin{array}{lllll}x_{1}&x_{2}&y_1&y_2\\
2&3&7&4&\\ 2&3&7&3&\\4&2&6&.5&\end{array} \right|=0$\\\\\\ that yields $-4x_1-16x_2+8y_1=0$.\\

The other weak defining hyperplanes are as follows:\\\\
$H_{3}$ is constructed on  $D9=\{1, 4, 5\}$ that yields $-7x_1+2y_1=0$.\\
$H_{4}$ is constructed on  $D10=\{2, 3, 6\}$, $D11=\{2, 3, 9\}$ and $D12=\{3, 6, 9\}$ and $D13=\{2, 6, 9\}$,\\ that yields $-27x_2+8y_1+6y_2=0$.\\
$H_{5}$ is constructed on $D14=\{3, 10, 9\}$ that yields $-6x_2+2y_1=0$.\\
$H_{6}$ is constructed on $D15=\{2, 6, 8\}$ that yields $-5x_2+2y_2=0$.\\
$H_{7}$ is constructed on $D14=\{2, 7, 8\}$ that yields $-5x_1-y_2=0$.\\

One can easily verify that the conditions of Theorem 4 are held and  $H_i, i=1, ..., 7$ are defining. It is worthwhile to note that the aforementioned PPS has only one strong defining hyperplane\\$H_8: -x_1-55x_2+17y_1+12y_2=0$; which is constructed on $D'=\{1, 2, 3\}$.\\

The following results can be attained by the aforementioned example:\\

\begin{itemize}
  \item There are one strong defining hyperplane and seven weak defining hyperplanes.
  \item The extreme ray binding from DMU $U_1$ is the intersection of four defining hyperplanes $H_1$, $H_2$, $H_3$ and $H_8$.
  \item The extreme ray binding from DMU $U_2$ is the intersection of five defining hyperplanes $H_1$, $H_4$, $H_6$, $H_7$ and $H_8$.
  \item The extreme ray binding from DMU $U_3$ is the intersection of four defining hyperplanes $H_2$, $H_4$, $H_5$ and $H_8$.
  \item All extreme rays are the intersection of strong and weak defining hyperplanes. Also all CCR-efficiency DMUs $U_1$, $U_2$ and $U_3$ lie on the weak defining hyperplanes.
\end{itemize}
\begin{table}
\caption{{\scriptsize Example 2. The results of evaluation CCR-efficient DMUs by models (\ref{tak}) and (\ref{tao}).}}
\begin{center}
\begin{tabular}{lccccccccccccccccccc}
  % after \\ : \hline or \cline{col1-col2} \cline{col3-col4} ...
\hline
DMU&&&&&&&l&&&&&&&&&q\\ \cline{6-9}  \cline{16-19}
&&&&&1&&&2&&&&&&&1&&&2\\ \hline
$D_1$&&&&&{\scriptsize FES}&&&{\scriptsize INFES}&&&&&&&{\scriptsize FES}&&&{\scriptsize INFES}& \\
$D_2$&&&&&{\scriptsize INFES}&&&{\scriptsize INFES}&&&&&&&{\scriptsize INFES}&&&{\scriptsize FES}&\\
$D_4$ &&&&&{\scriptsize INFES}&&&{\scriptsize FES}&&&&&&&{\scriptsize FES}&&&{\scriptsize INFES}& \\ \hline
\end{tabular}
\end{center}
\end{table}
\begin{table}
\caption{{\scriptsize Example 2. Weak efficient virtual DMUs. }}
\begin{center}
{\scriptsize
\begin{tabular}{lllllllllllllllllllllllllllllll}
\hline
DMU&&&&$D_{1}^{2}$&&&&$D_{1}^{4}$&&&&$D_{2}^{1}$&&&&$D_{2}^{2}$&&&&$D_{2}^{3}$&&&&$D_{4}^{1}$&&&&$D_{4}^{4}$\\
\hline
$x_1$&&&&2&&&&2&&&&2&&&&1&&&&1&&&&5&&&&4\\
$x_2$&&&&4&&&&3&&&&2&&&&3&&&&2&&&&2&&&&2\\
$y_1$&&&&7&&&&7&&&&3&&&&3&&&&2&&&&6&&&&6\\
$y_2$&&&&4&&&&3&&&&5&&&&5&&&&5&&&&1&&&&0\\
\hline
\end{tabular}}
\end{center}
\end{table}
\textit{4.3. Real word data}\\

We evaluated the data of 20 branchs of a bank in Iran using the proposed method. The data was previously analyzed by Amirteimoori et al. (2005), (see Table (\ref{ta})). Running the DEA model (\ref{CCRi}) (or (\ref{CCRo})) resulted in seven CCR-efficient units as 1, 4, 7, 12, 15, 17 and 20. Using the proposed method, the equations of weak defining hyperplanes were obtained as summarized below:\\\\
\begin{table}[t]
\caption{{\scriptsize Example 3. DMUs' data (extracted from [Amirteimoori et al. (2005), p. 689]).}}
\begin{center}\label{ta}
{\scriptsize
\begin{tabular}{llllllllllllllllllllllllllllll}
\hline
&&&input&&&&&&&&output\\ \cline{4-9}  \cline{12-18}
Branch&&&Staff&&&Computer &&Space m2&&&Deposits&&&Loans&&&Charge\\
&&&&&&terminals&&&&&&&&&&&\\
\hline
1&&&0.9503&&&0.70&&0.1550&&&0.1900&&&0.5214&&&0.2926\\
2&&&0.7962&&&0.60&&1.0000&&&0.2266&&&0.6274&&&0.4624\\
3&&&0.7982&&&0.75&&0.5125&&&0.2283&&&0.9703&&&0.2606\\
4&&&0.8651&&&0.55&&0.2100&&&0.1927&&&0.6324&&&1.0000\\
5&&&0.8151&&&0.85&&0.2675&&&0.2333&&&0.7221&&&0.2463\\
6&&&0.8416&&&0.65&&0.5000&&&0.2069&&&0.6025&&&0.5689\\
7&&&0.7189&&&0.60&&0.3500&&&0.1824&&&0.9000&&&0.7158\\
8&&&0.7853&&&0.75&&0.1200&&&0.1250&&&0.2340&&&0.2977\\
9&&&0.4756&&&0.60&&0.1350&&&0.0801&&&0.3643&&&0.2439\\
10&&&0.6782&&&0.55&&0.5100&&&0.0818&&&0.1835&&&0.0486\\
11&&&0.7112&&&1.00&&0.3050&&&0.2117&&&0.3179&&&0.4031\\
12&&&0.8113&&&0.65&&0.2550&&&0.1227&&&0.9225&&&0.6279\\
13&&&0.6586&&&0.85&&0.3400&&&0.1755&&&0.6452&&&0.2605\\
14&&&0.9763&&&0.80&&0.5400&&&0.1443&&&0.5143&&&0.2433\\
15&&&0.6845&&&0.95&&0.4500&&&1.0000&&&0.2617&&&0.0982\\
16&&&0.6127&&&0.90&&0.5250&&&0.1151&&&0.4021&&&0.4641\\
17&&&1.0000&&&0.60&&0.2050&&&0.0900&&&1.0000&&&0.1614\\
18&&&0.6337&&&0.65&&0.2350&&&0.0591&&&0.3492&&&0.0678\\
19&&&0.3715&&&0.70&&0.2375&&&0.0385&&&0.1898&&&0.1112\\
20&&&0.5827&&&0.55&&0.5000&&&0.1101&&&0.6145&&&0.7643\\
\hline
\end{tabular}}
\end{center}
\end{table}
The equations of defining hyperplanes binding from $DMU_1$:\\
\begin{enumerate}
  \tiny\item -972780000$x_1$  - 35830643110$x_3$+ 14965298195$y_1$+ 6971185000$y_2$=0
  \item - 8755020000$x_1$ - 322475787990$x_3$+ 134687683755$y_1$+ 62740665000$y_2$=0
  \item - 194556000$x_1$ - 7166128622$x_3$+ 2993059639$y_1$+ 1394237000$y_2$=0
  \item  -1359468485616$x_1$  - 28494222363367.38$x_3$+ 12093263476274$y_1$+ 5939111853909$y_2$+ 1073596674923.5$y_3$=0
  \item - 97278$x_2$ - 2645864$x_3$+ 1152988$y_1$+ 497000$y_2$=0
  \item - 875502$x_2$ - 23812776$x_3$+ 10376892$y_1$+ 4473000$y_2$=0
  \item - 194556$x_2$ - 52917280$x_3$+ 23059760$y_1$+ 9940000$y_2$=0
  \item - 1359468485616$x_2$- 22373297128800$x_3$+ 10143985960000$y_1$+ 4371989952000$y_2$+ 726507524000$y_3$=0
  \item - 943354$x_3$+ 388133$y_1$+ 139000$y_2$=0
  \item - 42450930$x_3$+ 17465985$y_1$+ 6255000$y_2$=0
  \item - 303325252038$x_3$+ 125752259800$y_1$+ 34411902300$y_2$+ 17703755450$y_3$=0
  \item - 286148$x_3$+ 96226$y_1$+ 50000$y_2$=0
  \item - 143074$x_3$+ 48113$y_1$+ 25000$y_2$=0
  \item - 100166891292$x_3$+ 26451830000$y_1$+ 17352331800$y_2$+ 49641649000$y_3$=0
\end{enumerate}
The equations of defining hyperplanes binding from $DMU_4$:\\
\begin{enumerate}
\tiny\item -135946848561600$x_1$-2849422236336738$x_3$+1209326347627400$y_1$+593911185390900$y_2$ + 107359667492350$y_3$=0
\item -303325252038$x_3$+ 125752259800$y_1$+3441190230$y_2$+17703755450$y_3$=0
\item -962972121000$x_1$-328314063906$x_2$+888323202000$y_1$+842460036000$y_3$=0
\item -69807437340$x_1$-109438021302$x_3$+474497788800$y_1$+535450401000$y_3$=0
\item -30930441000$x_1$-48184694626$x_2$+29114103500$y_2$+34847747500$y_3$=0
\item -1284492135536300$x_1$-338476850586989$x_3$+861567778754200$y_1$+351625903900200$y_2$ + 793901952483300$y_3$=0
\item -6201072000000$x_1$-63394549535400$x_3$+18021906285000$y_2$+7280349255000$y_3$=0
\item -3082802028$x_2$-1303063920$x_3$+2003476500$y_2$+702186000$y_3$=0
\item -9415572931881$x_2$-5266914542760$x_3$+9739254084500$y_1$+5725497127500$y_2$ + 787058521000$y_3$=0
\item -315203434864350$x_1$-1925392176272898$x_3$+911363930110150$y_1$+609134934067650$y_2$ + 116178086881850$y_3$=0
\item -402094816405800$x_1$-774926386660479$x_3$+506520288276700$y_1$+385152647434200$y_2$ + 169409773083050$y_3$=0
\item -469929571500$x_1$- 5215053344961$x_2$-3555135946650$x_3$+4096877084250$y_2$ + 1430528892750$y_3$=0
\item -49221399250365$x_1$- 102872514916490.925$x_2$-209914601960477.7$x_3$+188771599172697.5$y_1$ +133350621103897.5$y_2$+ 22536162160577.5$y_3$=0
\item -1396148746800$x_1$-218876042604$x_3$+948995577600$y_1$+1070900802000$y_3$=0
\item -498267154301400$x_1$-675598525370802$x_3$+412589517921900$y_2$+312004994380200$y_3$=0
\item -7193204732$x_2$-3040482480$x_3$+4674778500$y_2$+1638434000$y_3$=0
\item -203698200$x_1$-47097558$x_3$+186109800$y_3$=0
\item -209422312020$x_1$-328314063906$x_3$+142349336640$y_1$+160635120300$y_3$=0
\item -962972121000$x_1$-328314063906$x_2$+888323202000$y_1$+842460036000$y_3$=0
\item -2094223120200$x_1$-328314063906$x_3$+1423493366400$y_1$+1606351203000$y_3$=0
\item -588646116$x_3$+257626800$y_1$+73971000$y_3$=0
\item -203698200$x_1$-47097558$x_3$+186109800$y_3$=0
\item -20670240000$x_1$-211315165118$x_3$+60073020950$y_2$+24267830850$y_3$=0
\item -30930441000$x_1$-48184694626$x_2$+29114103500$y_2$+34847747500$y_3$=0
\item -156643190500$x_1$-1738351114987$x_2$-1185045315550$x_3$+1365625694750$y_2$ + 476842964250$y_3$=0
\item -49831698600$x_1$-67566609198$x_3$+41263078100$y_2$+31203619800$y_3$=0
\item -30930441000$x_1$-48184694626$x_2$+29114103500$y_2$+34847747500$y_3$=0
\item -125676913000$x_1$-48184694626$x_3$+49508977400$y_2$+87532406000$y_3$=0
\item -125676913000$x_1$-48184694626$x_3$+49508977400$y_2$+87532406000$y_3$=0
\item -1027600676$x_2$-434354640$x_3$+667825500$y_2$+234062000$y_3$=0
\item -44732808$x_2$+20631000$y_2$+11556000$y_3$=0
\item -203698200$x_1$-47097558$x_3$+186109800$y_3$=0
\item -203698200$x_1$-47097558$x_3$+186109800$y_3$=0
\item -320990707000$x_1$-109438021302$x_2$+296107734000$y_1$+280820012000$y_3$=0
\item -30930441000$x_1$-48184694626$x_2$+29114103500$y_2$+34847747500$y_3$=0
\item -30930441000$x_1$-48184694626$x_2$+29114103500$y_2$+34847747500$y_3$=0
\item -477909287649000$x_1$-338476850586989$x_2$+559061370687500$y_1$+207010083861500$y_2$ + 360957289402500$y_3$=0
\item -429988020954$x_2$+359622357000$y_1$+162565467000$y_2$+64387782000$y_3$=0
\item -777810000$x_1$-47097558$x_2$+931920000$y_3$=0
\item -30930441000$x_1$-48184694626$x_2$+29114103500$y_2$+34847747500$y_3$=0
\item -30930441000$x_1$-48184694626$x_2$+29114103500$y_2$+34847747500$y_3$=0
\item -30930441000$x_1$-48184694626$x_2$+29114103500$y_2$+34847747500$y_3$=0
\item -30930441000$x_1$-48184694626$x_2$+29114103500$y_2$+34847747500$y_3$=0
\item -203698200$x_1$-47097558$x_3$+186109800$y_3$=0
\item -2094223120200$x_1$-328314063906$x_3$+1423493366400$y_1$+1606351203000$y_3$=0
\item -196215372$x_3$+8587560+24657000$y_3$=0
\item -62010720000$x_1$-633945495354$x_3$+180219062850$y_2$+72803492550$y_3$=0
\item -31520343486435$x_1$-192539217627289.8$x_3$+91136393011015$y_1$+60913493406765$y_2$+11617808688185$y_3$=0
\item -49831698600$x_1$-67566609198$x_3$+41263078100$y_2$+31203619800$y_3$=0
\item -125676913000$x_1$-48184694626$x_3$+49508977400$y_2$+87532406000$y_3$=0
\item -203698200$x_1$-47097558$x_3$+186109800$y_3$=0
\item -588646116$x_2$+537594000$y_1$+220161000$y_3$=0
\item -89793064$x_3$+17110600$y_2$+8035800$y_3$=0
\item -196215372$x_2$+179198000$y_1$+73387000$y_3$=0
\item -269379192$x_3$+51331800$y_2$+24107400$y_3$=0
\item -600$x_3$+126$y_3$=0
\item -135946848561600$x_1$-2849422236336738$x_3$+1209326347627400$y_1$+593911185390900$y_2$+107359667492350$y_3$=0
\item -100166891292$x_3$+26451830000$y_1$+17352331800$y_2$+4964164900$y_3$=0
\item -1359468485616$x_2$-22373297128800$x_3$+10143985960000$y_1$+4371989952000$y_2$+726507524000$y_3$=0
\item -303325252038$x_3$+125752259800$y_1$+34411902300$y_2$+17703755450$y_3$=0
\end{enumerate}
The equations of defining hyperplanes binding from $DMU_7$:\\
\begin{enumerate}

 \tiny\item- 1284492135536300$x_1$       -338476850586989$x_3$ + 861567778754200$y_1$      + 351625903900200$y_2$+  793901952483300$y_3$=0
 \item- 18000$x_1$   - 54330$x_2$   + 50598$y_2$=0
 \item- 8890645750000$x_1$     - 4649133826135$x_2$ -16924059161800$x_3$     + 14494910015000$y_1$ + 13845016272500$y_2$=0
 \item- 2163150000$x_1$     - 1164847275$x_2$ -1944747000$x_3$      + 3260731500$y_2$=0
 \item- 62234520250000$x_1$     - 32543936782945$x_2$ -118468414132600$x_3$    + 101464370105000$y_1$ + 96915113907500$y_2$=0
 \item- 721050000$x_1$      - 388282425$x_2$ -648240000$x_3$      + 1086910500$y_2$=0
 \item- 4200$x_1$        - 12677$x_2$           + 11806$y_2$=0
 \item- 9415572931881$x_2$       -5266914542760$x_3$ + 9739254084500$y_1$       + 5725497127500$y_2$+  787058521000$y_3$=0
 \item - 3082802028$x_2$     -1303063920$x_3$  + 2003476500$y_2$+      702186000$y_3$=0
 \item- 800158117500$x_1$    - 4184220443521500$x_2$ -1523165324562$x_3$    + 1304541901350$y_1$ + 1246051464525$y_2$=0
 \item- 614548655100$x_1$       -818367903423$x_3$ + 609498131850$y_1$      + 685617039000$y_2$=0
 \item- 477982287300$x_1$      -636508369329$x_3$ + 474054102550$y_1$      + 533257697000$y_2$=0
 \item- 3605250000$x_1$       - 1941412125$x_2$ -3241245000$x_3$        + 5434552500$y_2$=0
 \item- 533438745000$x_1$     - 2789480295681000$x_2$ -1015443549708$x_3$     + 869694600900$y_1$ + 8307009763500000$y_2$=0
 \item- 614548655100$x_1$    -818367903423$x_3$ + 609498131850$y_1$ + 685617039000$y_2$=0
 \item- 409699103400$x_1$          -545578602282$x_3$ + 406332087900$y_1$      + 457078026000$y_2$=0
 \item- 402094816405800$x_1$      -774926386660479$x_3$ + 506520288276700$y_1$     + 385152647434200$y_2$+  169409773083050$y_3$=0
 \item- 1214119364$x_2$       -218213520$x_3$ + 1074571000$y_1$       + 676494000$y_2$=0
 \item- 9614000$x_1$       - 5177099$x_2$ -8643320$x_3$       + 14492140$y_2$=0
 \item- 40969295791344900$x_1$      -54557041860296577$x_3$ + 40632599291868150$y_1$      + 45707116982961000$y_2$=0
 \item- 8498835548$x_2$      -1527494640$x_3$ + 7521997000$y_1$      + 47354580000$y_2$=0
 \item- 596586144$x_2$    + 488586000$y_1$ + 298704000$y_2$=0
 \item- 25236750000$x_1$      - 13589884875$x_2$ -22688715000$x_3$       + 38041867500$y_2$=0
\item- 477982287300$x_1$       -636508369329$x_3$ + 474054102550$y_1$     + 533257697000$y_2$=0
 \item- 65362500$x_1$      -46889325$x_3$  + 70444850$y_2$=0
 \item- 24801090$x_1$      + 13235861$y_2$+  8266595$y_3$=0
 \item-57054$x_2$     + 33264$y_2$+  6000$y_3$=0
 \item- 1201750000$x_1$     - 647137375$x_2$ -1080415000$x_3$      + 1811517500$y_2$=0
 \item- 5602500$x_1$      -4019085$x_3$  + 6038130$y_2$=0
 \item- 8412250000$x_1$      - 4529961625$x_2$ -7562905000$x_3$      + 12680622500$y_2$=0
 \item- 171162$x_2$     + 99792$y_2$+  18000$y_3$=0
 \item- 233428299102$x_1$     + 127007416523$y_1$ + 93543291008$y_2$+    84460030097$y_3$=0
 \item- 767039328$x_2$   + 628182000$y_1$ + 384048000$y_2$=0
 \item- 12600$x_1$      - 38031$x_2$       + 3541860$y_2$=0
 \item- 6000$x_1$        - 1811000$x_2$           + 16866$y_2$=0
 \item- 409699103400$x_1$       -545578602282$x_3$ + 406332087900$y_1$      + 457078026000$y_2$=0
 \item- 5602500$x_1$    -4019085$x_3$  + 6038130$y_2$=0
 \item- 596586144$x_1$     + 299539709$y_1$ + 415833040$y_2$=0
 \item- 65362500$x_1$    -46889325$x_3$  + 70444850$y_2$=0
 \item- 767039328$x_1$     + 385122483$y_1$ + 534642480$y_2$=0
 \item - 10927074276$x_2$     -1963921680$x_3$ + 9671139000$y_1$    + 6088446000$y_2$=0
 \item - 70980$x_2$    + 42000$y_1$ + 38808$y_2$=0
 \item- 10140$x_2$    + 6000$y_1$ + 5544$y_2$=0
 \item- 30930441000$x_1$       - 48184694626$x_2$        + 291141035000$y_2$+  34847747500$y_3$=0
 \item- 156643190500$x_1$     - 1738351114987$x_2$ -1185045315550$x_3$      + 1365625694750$y_2$+  476842964250$y_3$=0
 \item- 49831698600$x_1$         -67566609198$x_3$  + 41263078100$y_2$+         31203619800$y_3$=0
 \item- 125676913000$x_1$      -48184694626$x_3$  + 49508977400$y_2$+      87532406000$y_3$=0
 \item - 1027600676$x_2$     -434354640$x_3$  + 667825500$y_2$+      234062000$y_3$=0
 \item - 44732808$x_2$     + 20631000$y_2$+  11556000$y_3$=0
 \item- 477909287649000$x_1$    - 338476850586989$x_2$  + 559061370687500$y_1$    + 207010083861500$y_2$+  360957289402500$y_3$=0
 \item - 429988020954$x_2$    + 359622357000$y_1$ + 162565467000$y_2$+    64387782000$y_3$=0
 \item- 7193204732$x_2$     -3040482480$x_3$  + 4674778500$y_2$+      1638434000$y_3$=0
 \item- 469929571500$x_1$     - 5215053344961$x_2$ -3555135946650$x_3$      + 4096877084250$y_2$+  1430528892750$y_3$=0
 \item- 492213992503650$x_1$- 1028725149164909.25$x_2$ -2099146019604777$x_3$+ 1887715991726975$y_1$ + 1333506211038975$y_2$+  225361621605775$y_3$=0
\end{enumerate}
The equations of defining hyperplanes binding from $DMU_{12}$:\\
\begin{enumerate}
\tiny\item- 4922139925036500$x_1$-10287251491649092.5$x_2$-20991460196047770$x_3$+18877159917269750$y_1$+13335062110389750$y_2$ +2253616216057750$y_3$=0
\item- 20670240000000$x_1$- 211315165118000$x_3$+60073020950000$y_2$+24267830850000$y_3$=0
\item- 49831698600000$x_1$-67566609198000$x_3$+41263078100000$y_2$+31203619800000$y_3$=0
\item- 15664319050000000$x_1$- 173835111498700000$x_2$-118504531555000000$x_3$+136562569475000000$y_2$ +47684296425000000$y_3$=0
\item- 315203434864350$x_1$-1925392176272898$x_3$+911363930110150$y_1$ +609134934067650$y_2$ +116178086881850$y_3$=0
\item- 6201072000000$x_1$-63394549535400$x_3$+18021906285000$y_2$+7280349255000$y_3$=0
\item- 4699295715000000$x_1$-52150533449610000$x_2$- 35551359466500000$x_3$+40968770842500000$y_2$+14305288927500000$y_3$=0
\item- 84122500000000000$x_1$-45299616250000000$x_2$- 75629050000000000$x_3$+126806225000000000$y_2$=0
\item- 56025000000000$x_1$- 40190850000000$x_3$+60381300000000$y_2$ =0
\item- 65362500000000$x_1$-46889325000000$x_3$+70444850000000$y_2$ =0
\item- 477982287300000$x_1$-636508369329000$x_3$+474054102550000$y_1$ +533257697000000$y_2$=0
\item- 622345202500000000$x_1$-325439367829450000$x_2$-1184684141326000000$x_3$+1014643701050000000$y_1$ + 969151139075000000$y_2$=0
\item- 25236750000000000$x_1$-13589884875000000$x_2$-22688715000000000$x_3$+ 38041867500000000$y_2$=0
\item- 61454865510000$x_1$-81836790342300$x_3$+60949813185000$y_1$+68561703900000$y_2$=0
\item- 72105000000000000$x_1$-38828242500000000$x_2$-64824900000000000$x_3$+108691050000000000$y_2$=0
\item- 21631500000000000$x_1$-11648472750000000$x_2$-19447470000000000$x_3$+ 2607315000000000$y_2$=0
\item- 80015811750000000$x_1$-41842204435215000$x_2$-152316532456200000$x_3$+130454190135000000$y_1$+124605146452500000$y_2$=0
\item- 88906457500000000$x_1$- 46491338261350000$x_2$-169240591618000000$x_3$+ 144949100150000000$y_1$+138450162725000000$y_2$=0
\item- 12017500000000000$x_1$-6471373750000000$x_2$-10804150000000000$x_3$+ 18115175000000000$y_2$=0
\item- 533438745000000000$x_1$-278948029568100000$x_2$-1015443549708000000$x_3$+ 869694600900000000$y_1$+830700976350000000$y_2$=0
\item- 3605250000000000$x_1$-1941412125000000$x_2$-3241245000000000$x_3$+5434552500000000$y_2$ =0
\item- 144691680000000000 $x_1$-1479206155826000000$x_3$+420511146650000000$y_2$ +169874815950000000$y_3$=0
\item- 402094816405800$x_1$-774926386660479$x_3$+506520288276700$y_1$+385152647434200$y_2$+169409773083050$y_3$=0
\item- 4165274217900000000$x_1$-5546715789867000000$x_3$+4131042893650000000$y_1$+ 4646959931000000000$y_2$=0
\item- 13219250000000000$x_1$-7118511125000000$x_2$-11884565000000000$x_3$+19926692500000000$y_2$=0
\item- 409699103400000$x_1$- 545578602282000$x_3$+406332087900000$y_1$+457078026000000$y_2$=0
\item- 39532500000000$x_1$- 66720000000000$x_3$+53210100000000$y_2$=0
\item- 291662435700000$x_1$-752132978400000$x_3$+431589130800000$y_1$+407006674500000$y_2$=0
\item- 11859750000000$x_1$-20016000000000$x_3$+ 15963030000000$y_2$ =0
\item- 43749365355000$x_1$-112819946760000$x_3$+64738369620000$y_1$+61051001175000$y_2$=0
\item- 48610405950000$x_1$-125355496400000$x_3$+71931521800000$y_1$+67834445750000$y_2$=0
\item- 1976625000000$x_1$-3336000000000$x_3$+ 2660505000000$y_2$=0
\item- 6588750000000$x_1$-11120000000000$x_3$+8868350000000$y_2$=0
\item- 19766250000000000$x_1$-33360000000000000$x_3$+26605050000000000$y_2$=0
\end{enumerate}
The equations of defining hyperplanes binding from $DMU_{15}$:\\
\begin{enumerate}
\tiny\item- 150697636$x_1$+93188442$y_1$+101467310$y_3$=0
 \item- 452092908$x_1$+279565326$y_1$+304401930$y_3$=0
 \item- 18000$x_1$+12321$y_1$=0
 \item- 180$x_3$+81$y_1$=0
 \item- 291657574659405$x_1$-752120442850360$x_3$ + 431581937647820$y_1$   + 406999891055425$y_2$ =0
 \item- 43749365355$x_1$-112819946760$x_3$ + 64738369620$y_1$ + 61051001175$y_2$ =0
 \item - 10927074276$x_2$-1963921680$x_3$+9671139000$y_1$+6088446000$y_2$ =0
 \item- 4861040595$x_1$-12535549640$x_3$+7193152180$y_1$+6783444575$y_2$  =0
 \item- 2916624357$x_1$-7521329784$x_3$+4315891308$y_1$+4070066745$y_2$=0
 \item- 5334387450000$x_1$-2789480295681$x_2$-10154435497080$x_3$+8696946009000$y_1$+8307009763500$y_2$ =0
 \item- 409699103400$x_1$-545578602282$x_3$+406332087900$y_1$+457078026000$y_2$=0
 \item- 8890645750000$x_1$-4649133826135$x_2$-16924059161800$x_3$+14494910015000$y_1$+13845016272500$y_2$ =0
 \item- 233428299102$x_1$+127007416523$y_1$+93543291008$y_2$+84460030097$y_3$=0
 \item- 767039328$x_1$+385122483$y_1$+534642480$y_2$=0
 \item - 1214119364$x_2$-218213520$x_3$+1074571000$y_1$+676494000$y_2$ =0
 \item- 80015811750000$x_1$-41842204435215$x_2$-152316532456200$x_3$+130454190135000$y_1$+124605146452500$y_2$ =0
 \item- 614548655100$x_1$-818367903423$x_3$+609498131850$y_1$+685617039000$y_2$=0
 \item - 8498835548$x_2$-1527494640$x_3$+7521997000$y_1$+4735458000$y_2$ =0
 \item- 596586144$x_1$+299539709$y_1$+415833040$y_2$ =0
 \item- 62234520250000$x_1$-32543936782945$x_2$-118468414132600$x_3$+101464370105000$y_1$+96915113907500$y_2$ =0
 \item- 477982287300$x_1$-636508369329$x_3$+474054102550$y_1$+533257697000$y_2$ =0
 \item - 767039328$x_2$+628182000$y_1$+384048000$y_2$ =0
 \item - 596586144$x_2$+488586000$y_1$+298704000$y_2$ =0
 \item- 402094816405800$x_1$-774926386660479$x_3$+506520288276700$y_1$+385152647434200$y_2$+169409773083050$y_3$=0
 \item -49221399250365$x_1$- 102872514916490.925$x_2$-209914601960477.7$x_3$+188771599172697.5$y_1$ +133350621103897.5$y_2$+ 22536162160577.5$y_3$=0
 \item- 477909287649000$x_1$-338476850586989$x_2$+559061370687500$y_1$+207010083861500$y_2$+360957289402500$y_3$=0
 \item - 9415572931881$x_2$-5266914542760$x_3$+9739254084500$y_1$+5725497127500$y_2$+787058521000$y_3$=0
 \item- 320990707000$x_1$-109438021302$x_2$+296107734000$y_1$+280820012000$y_3$=0
 \item- 962972121000$x_1$-328314063906$x_2$+888323202000$y_1$+842460036000$y_3$=0
 \item- 698074373400$x_1$-109438021302$x_3$+474497788800$y_1$+535450401000$y_3$=0
 \item - 196215372$x_2$+179198000$y_1$+73387000$y_3$=0
 \item- 1284492135536300$x_1$-338476850586989$x_3$+861567778754200$y_1$+351625903900200$y_2$+793901952483300$y_3$=0
 \item - 429988020954$x_2$+359622357000$y_1$+162565467000$y_2$+64387782000$y_3$=0
 \item - 588646116$x_2$+537594000$y_1$+220161000$y_3$=0
 \item- 196215372$x_3$+85875600$y_1$+24657000$y_3$=0
 \item- 315203434864350$x_1$-1925392176272898$x_3$+911363930110150$y_1$+609134934067650$y_2$+116178086881850$y_3$=0
 \item- 2094223120200$x_1$-328314063906$x_3$+1423493366400$y_1$+1606351203000$y_3$=0
 \item - 588646116$x_3$+257626800$y_1$+73971000=0
 \item- 135946848561600$x_1$-2849422236336738$x_3$+1209326347627400$y_1$+593911185390900$y_2$+107359667492350$y_3$=0
 \item - 194556$x_2$-5291728$x_3$+2305976$y_1$+994000$y_2$ =0
 \item - 97278$x_2$-2645864$x_3$+1152988$y_1$+497000$y_2$ =0
 \item- 972780000$x_1$-35830643110$x_3$+14965298195$y_1$+6971185000$y_2$ =0
 \item - 875502$x_2$-23812776$x_3$+10376892$y_1$+4473000$y_2$ =0
 \item - 42450930$x_3$+17465985$y_1$+6255000$y_2$ =0
 \item - 1359468485616$x_2$-22373297128800$x_3$+10143985960000$y_1$+4371989952000$y_2$+726507524000$y_3$=0
 \item - 303325252038$x_3$+125752259800$y_1$+34411902300$y_2$+17703755450$y_3$=0
 \item- 194556000$x_1$-7166128622$x_3$+2993059639$y_1$+1394237000$y_2$ =0
 \item - 943354$x_3$+388133$y_1$+139000$y_2$ =0
\end{enumerate}
The equations of defining hyperplanes binding from $DMU_{17}$:\\
\begin{enumerate}
\tiny\item- 972780000$x_1$-35830643110$x_3$+14965298195$y_1$+6971185000$y_2$ =0
 \item- 8890645750000$x_1$+4649133826135$x_2$-16924059161800$x_3$+14494910015000$y_1$+13845016272500$y_2$ =0
 \item- 8755020000$x_1$-322475787990$x_3$+134687683755$y_1$+62740665000$y_2$ =0
 \item- 80015811750000$x_1$+41842204435215$x_2$-152316532456200$x_3$+130454190135000$y_1$+124605146452500$y_2$ =0
 \item- 721050000$x_1$+388282425$x_2$-648249000$x_3$+1086910500$y_2$ =0
\item- 658875$x_1$-1112000$x_3$+886835$y_2$=0
 \item- 62234520250000$x_1$+32543936782945$x_2$-118468414132600$x_3$+101464370105000$y_1$+96915113907500$y_2$ =0
 \item- 62010720000$x_1$-633945495354$x_3$+180219062850$y_2$+72803492550$y_3$=0
\item- 6000$x_1$+18110$x_2$+16866$y_2$ =0
 \item- 5334387450000$x_1$+2789480295681$x_2$-10154435497080$x_3$+8696946009000$y_1$+8307009763500$y_2$ =0
 \item- 492213992503650$x_1$+1028725149164909.25$x_2$-2099146019604777$x_3$+1887715991726975$y_1$+1333506211038975$y_2$+225361621605775$y_3$=0
\item- 4861040595$x_1$-12535549640$x_3$+7193152180$y_1$+6783444575$y_2$=0
\item- 469929571500$x_1$+5215053344961$x_2$-3555135946650$x_3$+4096877084250$y_2$+1430528892750$y_3$=0
 \item- 43749365355$x_1$-112819946760$x_3$+64738369620$y_1$+61051001175$y_2$=0
 \item- 42000$x_1$+126770$x_2$+118062$y_2$=0
\item- 395325$x_1$-667200$x_3$+532101$y_2$=0
 \item- 3605250000$x_1$+1941412125$x_2$-3241245000$x_3$+5434552500$y_2$=0
 \item- 315203434864350$x_1$-1925392176272898$x_3$+911363930110150$y_1$+609134934067650$y_2$+116178086881850$y_3$=0
\item- 2916624357$x_1$-7521329784$x_3$+4315891308$y_1$+4070066745$y_2$=0
 \item- 2916580607634645$x_1$-7521216964053240$x_3$+4315826569630380$y_1$+4070005693998825$y_2$ =0
 \item- 25236750000$x_1$+13589884875$x_2$-22688715000$x_3$+38041867500$y_2$=0
 \item- 2163150000$x_1$+1164847275$x_2$-1944747000$x_3$+3260731500$y_2$=0
 \item- 20670240000$x_1$-211315165118$x_3$+60073020950$y_2$+24267830850$y_3$=0
 \item- 206681729760000$x_1$-2112940336014882$x_3$+600670136479050$y_2$+242654040669150$y_3$=0
\item- 1976625$x_1$-3336000$x_3$+2660505$y_2$=0
 \item- 194556000$x_1$-7166128622$x_3$+2993059639$y_1$+1394237000$y_2$=0
 \item- 18000$x_1$+5433$x_2$+50598$y_2$=0
 \item- 135946848561600$x_1$-2849422236336738$x_3$+1209326347627400$y_1$+593911185390900$y_2$+107359667492350$y_3$=0
\item- 126000$x_1$+380310$x_2$+354186$y_2$=0
 \item- 1201750000$x_1$+647137375$x_2$-1080415000$x_3$+1811517500$y_2$ =0
 \item- 11859750000000$x_1$-20016000000000$x_3$+15963030000000$y_2$=0
 \item- 10700580000$x_1$-394137074210$x_3$+164618280145$y_1$+76683035000$y_2$ =0
 \item- 97278$x_2$-2645864$x_3$+1152988$y_1$+497000$y_2$ =0
\item- 9415572931881$x_2$-5266914542760$x_3$+9739254084500$y_1$+5725497127500$y_2$+787058521000$y_3$=0
\item- 875502$x_2$-23812776$x_3$+10376892$y_1$+4473000$y_2$ =0
 \item- 8498835548$x_2$-1527494640$x_3$+7521997000$y_1$+4735458000$y_2$ =0
 \item- 74061281204$x_2$-13311024720$x_3$+65548831000$y_1$+41266134000$y_2$ =0
 \item- 70980$x_2$+42000$y_1$+38808$y_2$ =0
 \item- 57054$x_2$+33264$y_2$+6000$y_3$=0
 \item- 3082802028$x_2$-1303063920$x_3$+2003476500$y_2$+702186000$y_3$=0
 \item- 300$x_2$+180$y_2$=0
 \item- 194556$x_2$-5291728$x_3$+2305976$y_1$+994000$y_2$ =0
 \item- 171162$x_2$+99792$y_2$+18000$y_3$=0
 \item- 1359468485616$x_2$-22373297128800$x_3$+10143985960000$y_1$+4371989952000$y_2$+726507524000$y_3$=0
 \item- 1214119364$x_2$-218213520$x_3$+1074571000$y_1$+676494000$y_2$ =0
 \item- 10927074276$x_2$-1963921680$x_3$+9671139000$y_1$+6088446000$y_2$ =0
 \item- 10700580$x_2$-291045040$x_3$+126828680$y_1$+54670000$y_2$ =0
 \item- 1027600676$x_2$-434354640$x_3$+667825500$y_2$+234062000$y_3$=0
 \item- 10140$x_2$+6000$y_1$+5544$y_2$ =0
 \item - 89793064$x_3$+17110600$y_2$+8035800$y_3$=0
 \item - 3000$x_3$+615$y_2$ =0
 \item - 2693791920$x_3$+51331800$y_2$+24107400$y_3$=0
 \item - 143074$x_3$+48113$y_1$+25000$y_2$ =0
 \item - 100166891292$x_3$+26451830000$y_1$+17352331800$y_2$+4964164900$y_3$=0
\end{enumerate}
The equations of defining hyperplanes binding from $DMU_{20}$:\\
\begin{enumerate}
\tiny\item- 1284492135536300$x_1$-338476850586989$x_3$+861567778754200$y_1$+351625903900200$y_2$+793901952483300$y_3$=0
 \item- 698074373400$x_1$-109438021302$x_3$+474497788800$y_1$+535450401000$y_3$=0
 \item- 125676913000$x_1$-48184694626$x_3$+49508977400$y_2$+87532406000$y_3$=0
 \item- 30930441000$x_1$-48184694626$x_2$+29114103500$y_2$+34847747500$y_3$=0
 \item- 77781000$x_1$-47097558$x_2$+93192000$y_3$=0
 \item- 962972121000$x_1$-328314063906$x_2$+888323202000$y_1$+842460036000$y_3$=0
 \item- 2094223120200$x_1$-328314063906$x_3$+1423493366400$y_1$+1606351203000$y_3$=0
 \item- 320990707000$x_1$-109438021302$x_2$+296107734000$y_1$+280820012000$y_3$=0
 \item- 477909287649000$x_1$-338476850586989$x_2$+559061370687500$y_1$+207010083861500$y_2$+360957289402500$y_3$=0
 \item- 248010900$x_1$+132358610$y_2$+82665950$y_3$=0
 \item- 233428299102$x_1$+127007416523$y_1$+93543291008$y_2$+84460030097$y_3$=0
 \item- 452092908$x_1$+279565326$y_1$+304401930$y_3$=0
 \item- 150697636$x_1$+93188442$y_1$+101467310$y_3$=0
 \item- 45858$x_1$+34962$y_3$=0
 \item- 203698200$x_1$-47097558$x_3$+186109800$y_3$=0
 \item- 1396148746800$x_1$-218876042604$x_3$+948995577600$y_1$+1070900802000$y_3$=0
\end{enumerate}
It is worthwhile to note that there are 216 weak defining hyperplanes and there is only one strong defining hyperplane that is constricted on DMUs 4, 7, 12, 15, 17.
\section{Conclusions}
Until now, less attention has been paid regarding finding {\it weak} facet of PPS of DEA models (see Wei et al. (2007)). Following Jahanshahloo et al. (2007), in this paper we proposed a method for finding all weak defining hyperplanes of the PPS of the CCR model. To do this, the performance of each DMUs was firstly evaluated using models (\ref{CCRi}) or (\ref{CCRo}), and all CCR-inefficient cases from the PPS were then removed. By introducing a variance of super-efficient models (see models (\ref{tak}) and (\ref{tao})) and using properties 2-7, the weak efficient virtual DMUs and the strong efficient DMUs are found. A supporting hyperplane was found to be a weak defining hyperplane if at least one weak efficient virtual DMU lies on it. Using the proposed method, one can check which CCR-efficient DMUs lie on the extreme rays (edges) of the PPS of the CCR model; which extreme DMUs lie on the weak defining hyperplanes, and how many defining hyperplanes they are on. In addition, these hyperplanes are useful in sensitivity and stability analysis. Our algorithm can easily be implemented using existing packages for operation research, such as GAMS.

%\section*{Acknowledgement}
%The authors are very grateful to the two anonymous referees and the editor, Prof. M. Cross for their comments
%and suggestions which have been very helpful in improving this paper.

\section*{References}
Amirteimoori A., Kordrostami S., 2005. ``Efficient surfaces and an efficiency index in DEA: a constant

returns to scale", \textit{Applied Mathematics and Computation }, 163 683-691.\\
Amirteimoori A., Kordrostami S., 2012. ``Generating strong defining hyperplanes of the production

possibility set in data envelopment analysis". \textit{Applied Mathematics Letters}, 25 (3) 605-609.\\
Banker R.D., 1984. ``Estimation most productive scale size using data envelopment analysis".\textit{European}

\textit{Journal of Operational Research}, 17 35-44.\\
Banker R.D., Thrall R.M., 1992. ``Estimation of returns to scale using data envelopment analysis".

\textit{European Journal of Operational Research}, 62 (2) 74-84.\\
Banker R.D., Charnes A., Cooper W.W., 1984. ``Some models for estimating technical and scale

inefficiencies in data envelopment analysis". \textit{Management Science} 30 (9) 1078-1092.\\
Bazaraa M., Shetty C.M., 1990. Linear Programming and Network Flows, John Wiley $\&$ Sons.\\
Charnes A., Cooper W.W., Rodes E., 1978. ``Measuring the efficiency of decision making units'', \textit{European}

\textit{Journal of Operational Research}, 2 (6) 429-444.\\
Cooper W.W., Seiford L.M., Tone K., 2000. ``Data Envelopment Analysis: a comprehensive text with

models, applications, references, and DEA-Solver software'', Kluwer Academic Publishers.\\
Davtalab-Olyaie M., Roshdi I., Partovi Nia V., Asgharian D., 2014. ``On characterizing full

dimensional weak facets in DEA with variable returns to scale technology". \textit{Optimization},

163 483-488.\\
Fukuyama H., Mirdehghan S.M., 2012. ``Identifying the efficiency status in network DEA''. \textit{European}

\textit{Journal of Operational Research} 220 85-92.\\
Hosseinzadeh Lotfi F., Jahanshahloo G.R., Mozaffari M.R., Gerami J., 2011. ``Finding DEA-efficient

hyperplanes using MOLP efficient faces". \textit{Journal of Computational and Applied Mathematics}, 235

1227-1231.\\
Jahanshahloo G.R., Hosseinzadeh Lotfi F., Shoja N., Sanei M., Tohidi G., 2005a. ``Sensitivity and

stability analysis in DEA". \textit{Applied Mathematics and Computation }, 169 897-904.\\
Jahanshahloo G.R., Hosseinzadeh Lotfi F., Zhiani Rezai H. and Rezai Balf F., 2007. ``Finding Strong

Defining Hyperplanes of Production Possibility Set". \textit{European Journal of Operational Research} 25(3)

605-609.\\
Jahanshahloo G.R., Hosseinzadeh Lotfi F., Zohrehbandian M., 2005b. ``Finding the piecewise linear

frontier production function in data envelopment analysis". \textit{Applied Mathematics and Computation},

163 483-488.\\
Jahanshahloo G.R., Shirzadi A., Mirdehghan S.M., 2009. ``Finding strong defining hyperplanes of PPS

using multiplier form". \textit{European Journal of Operational Research}, 194 (3) 933-938.\\
Jahanshahloo G.R., Hosseinzadeh Lotfi F., Akbarian D., 2010. ``Finding weak defining hyperplanes of

PPS of the BCC model". \textit{Applied Mathematical Modelling} 34 (11) 3321-3332.\\
Olesen O.B. and Petersen N.C., 2003. ``Identification and use of efficient faces and facets in DEA".

\textit{Journal of Productivity Analysis}, 20 323-360.\\
Seiford L.M. and Thrall R.M., 1990. ``Recent developments in DEA, The mathematical programming

approach to frontier analysis". \textit{Journal of Econometrics}. 46 7-38.\\
Wei Q.L., Yan H. and Hao G., 2007. ``Characteristics and structures of weak efficient surfaces of

production possibility sets". \textit{Journal of Mathematics Analysis and Application}. 327 1055-1074.\\
Yu G., Wei Q., Brockett P. and Zhou L., 1996. ``Construction of all DEA efficient surfaces of the

production possibility set under the Generalized Data Envelopment Analysis Model.

\textit{European Journal of Operational Research}, 95 (3) 491-510.
\newpage

\begin{figure}
[b]
\unitlength=4.5mm
\begin{center}
\begin{picture}(5,-33)(0,-2)
\put(0,0){\line(1,0){6}} \put(0,0){\line(0,1){4}}
\thicklines\put(0,4){\line(-1,-3){2}}
\thicklines\put(0,0){\line(-1,-1){2}}
\thicklines\put(6,0){\line(1,-1){2}}
\put(-2.2,-2.2){{\circle*{.4}}\makebox(-1,2)[]{\shortstack[b]{1}}}
\put(8.2,-2.2){{\circle*{.4}}\makebox(1,1)[]{\shortstack[b]{6}}}
\put(6,4){{\circle*{.4}}\makebox(0,1)[]{\shortstack[b]{5}}}
\put(0,0){{\circle*{.4}}\makebox(1,1)[]{\shortstack[b]{2}}}
\put(6,0){{\circle*{.4}}\makebox(-2,1)[]{\shortstack[b]{4}}}
\put(0,4){{\circle*{.4}}\makebox(-2,1)[]{\shortstack[b]{3}}}
\thicklines \put(0,0){\line(1,0){6}}\put(6,0){\line(0,1){4}}
\put(0,0){\line(0,1){4}} \put(6,4){\line(1,-3){2}}
\put(0,4){\line(1,0){6}}
\end{picture}
\end{center}
\caption{Example 1.}\label{4} \end{figure}
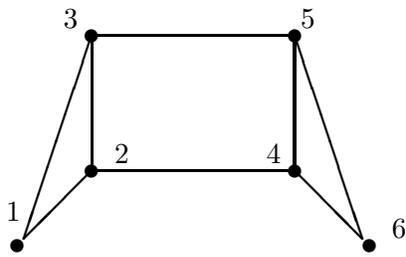\
\end{document}